\documentclass{styles/svproc}
\usepackage{url}

\usepackage[hidelinks]{hyperref}
\usepackage{type1cm}        % activate if the above 3 fonts are
\usepackage{makeidx}         % allows index generation
\usepackage{graphicx}        % standard LaTeX graphics tool
\usepackage{multicol}        % used for the two-column index
\usepackage[bottom]{footmisc}% places footnotes at page bottom

\usepackage{newtxtext}       %
\usepackage[varvw]{newtxmath}       % selects Times Roman as basic font

\usepackage{bm}% for bold math symbols
\usepackage{siunitx}% for table alignments

\usepackage{amsmath}
\usepackage{mathtools}
\usepackage{float}
\usepackage{subcaption}

\newcounter{algorithm}%[section]

\newcommand{\N}{{\mathbb{N}}} % natural numbers {1, 2, ...}
\newcommand{\R}{{\mathbb{R}}} % reals
\DeclareSymbolFont{bbold}{U}{bbold}{m}{n}
\DeclareSymbolFontAlphabet{\mathbbold}{bbold}

\usepackage{cite}

\newcommand{\dd}{\mathrm{d}}

\begin{document}

\mainmatter              % start of a contribution
\title{Low-rank variance reduction for uncertain radiative transfer with control variates}
\titlerunning{Low-rank for uncertain radiative transfer}  % abbreviated title (for running head)
\author{Chinmay Patwardhan\inst{1}\footnote{Equal contribution.} \and Pia Stammer\inst{2}$^\dagger$ \and
Emil L{\o}vbak\inst{1} \and Jonas Kusch\inst{3} \and Sebastian Krumscheid\inst{1}}
\authorrunning{C. Patwardhan, P. Stammer, et al.} % abbreviated author list (for running head)
\tocauthor{Chinmay Patwardhan, Pia Stammer, Emil L{\o}vbak, Jonas Kusch,
Sebastian Krumscheid}
\institute{Karlsruhe Institute of Technology (KIT), %Institute of Applied and Numerical Mathematics, 
Karlsruhe, Germany\\
\email{\{chinmay.patwardhan, emil.loevbak, sebastian.krumscheid\}@kit.edu}
\and
Delft University of Technology (TU Delft), %Department of Radiation Science \& Technology, 
Delft, Netherlands\\
\email{p.k.stammer@tudelft.nl}
\and
Norwegian University of Life Sciences (NMBU), %Scientific Computing, 
Ås, Norway\\ 
\email{jonas.kusch@nmbu.no}
}

\maketitle              % typeset the title of the contribution

\begin{abstract}
The radiative transfer equation models various physical processes ranging from plasma simulations to radiation therapy. In practice, these phenomena are often subject to uncertainties. Modeling and propagating these uncertainties requires accurate and efficient solvers for the radiative transfer equations. Due to the equation's high-dimensional phase space, fine-grid solutions of the radiative transfer equation are computationally expensive and memory-intensive. In recent years, dynamical low-rank approximation has become a popular method for solving kinetic equations due to the development of computationally inexpensive, memory-efficient and robust algorithms like the augmented basis update \& Galerkin integrator. In this work, we propose a low-rank Monte Carlo estimator and combine it with a control variate strategy based on multi-fidelity low-rank approximations for variance reduction. We investigate the error analytically and numerically and find that a joint approach to balance rank and
grid size is necessary. Numerical experiments further show that the efficiency of estimators can be improved using dynamical low-rank approximation, especially in the context of control variates.
\end{abstract}
%
%%%%%%%%%%%%%%%%%%%%%%%%%%%%%%%%%%%%%%%%%%%%%%%%%%%%%%%%%%%%%%%%%%%%%%%%%%
\section{Introduction}
%%%%%%%%%%%%%%%%%%%%%%%%%%%%%%%%%%%%%%%%%%%%%%%%%%%%%%%%%%%%%%%%%%%%%%%%%%

Kinetic equations play a key role in modeling various natural phenomena from plasma physics~\cite{howes_kinetic_2008,zweibel_magnetic_2009,einkemmer_semi-lagrangian_2023} to radiation therapy planning~\cite{duclous_deterministic_2010,kusch_robust_2023}. When performing numerical simulations, one often encounters uncertainties due to inaccurate measurements or unknown parameters. In forward uncertainty quantification (UQ), one assumes a probability distribution on these uncertainties and computes quantities of interest, e.g., the expectation or variance, of the corresponding computational results. Despite advances in both computing power and efficient solution methods, numerically solving such equations with uncertainties is still a computationally expensive and memory-intensive task. Hence, a combined approach of efficient sampling methods and efficient solvers is required. To this end, we present a novel combination of hierarchical Monte Carlo sampling methods with the dynamical low-rank approximation \cite{koch_dynamical_2007}. 

We consider the radiative transfer equation (RTE) which models the transport of particles through a background material
\begin{subequations}\label{eq:DLRAMC_3Dradiativetransport}
\begin{equation}
    \partial_{t} \psi(t, x, \Omega) + \Omega \cdot \nabla_{x} \psi(t, x, \Omega) = \sigma_s(x) \left( \frac{1}{4\pi}\int_{\mathbb{S}^2} \psi(t, x, \Omega^\prime) \dd \Omega^\prime - \psi(t, x, \Omega) \right),
    \end{equation}
    \begin{equation}
    \psi(t=t_{0},x,\Omega)= \psi_0(x,\Omega).            
    \end{equation}
\end{subequations}
Here, $\psi(t, x, \Omega)$ denotes the particle density (or angular flux) at time $t\in \mathbb{R}_{\geq 0}$, spatial position $x\in \mathcal{D}_x$ and travelling in direction $\Omega \in \mathbb{S}^2$. The particle density modeled by \eqref{eq:DLRAMC_3Dradiativetransport} is subject to discontinuous changes in direction due to collisions with the background material with scattering rate $\sigma_{s}(x)$. Moreover, we assume that the particle density does not reach the domain boundary. Due to the 6-dimensional phase space, numerically solving problems involving equations of the form \eqref{eq:DLRAMC_3Dradiativetransport} quickly becomes computationally expensive. Hence, a multitude of numerical methods have been developed to tackle such problems, e.g., particle-based Monte Carlo simulation~\cite{bird_monte_1978} and spectral methods, such as moment models~\cite{grad_kinetic_1949}. For an overview of methods, we refer to~\cite{dimarco_numerical_2014} and references therein.

Recently, dynamical low-rank approximation \cite{koch_dynamical_2007} has been widely used to solve kinetic equations in plasma physics \cite{einkemmer_low-rank_2018,coughlin_efficient_2022}, thermal radiative transfer \cite{baumann_energy_2024,frank_asymptotic-preserving_2025,patwardhan_parallel_2025} as well as problems related to nuclear power \cite{kusch_low-rank_2022} or medical physics \cite{kusch_robust_2023,stammer_deterministic_2024}. The core idea of dynamical low-rank approximation (DLRA) for radiative transfer is to restrict the underlying radiative transfer equation to low-rank solutions. This is achieved by projecting the dynamics onto the tangent space of the manifold of low-rank functions, leading to novel evolution equations for the low-rank factors of the solution. However, the rank is often overestimated to ensure it is sufficiently high. Overestimated ranks lead to ill-conditioned evolution equations \cite{koch_dynamical_2007}, requiring the derivation of novel time integrators tailored to the underlying structure of the manifold of low-rank solutions. The most frequently used robust time integrators for DLRA include projector--splitting \cite{lubich_projector-splitting_2014,kieri_discretized_2016} and basis-update \& Galerkin (BUG) integrators \cite{ceruti_unconventional_2022,ceruti_rank-adaptive_2022,ceruti_parallel_2024,ceruti_robust_2024,kusch_second-order_2024}. Although numerous works use dynamical low-rank approximation to simulate deterministic kinetic problems efficiently, approaches to include uncertainties as an additional dimension of the phase space have also been explored \cite{sapsis_dynamically_2009,musharbash_symplectic_2020,kusch_dynamical_2022,stammer_uncertainty_2023}. Here, most notably, \cite{stammer_uncertainty_2023} uses a DLRA tensor integrator to include uncertainties into kinetic simulations. However, such approaches require an intrusive modification of a solver's code. Further, especially in higher dimensions their efficiency depends strongly on the design of the tensor structures and is not yet well explored.

Monte Carlo methods are frequently used for forward uncertainty quantification due to their non-intrusive nature and ease of parallelization. They also tend to be favorable for higher-dimensional uncertainties, as they avoid the construction of parameter-space grids. However, direct Monte Carlo sampling is often prohibitively expensive in practical applications. Hence, variance reduction techniques such as control variates \cite{wilson_variance_1984} or Multilevel Monte Carlo methods (MLMC) \cite{margenov_multilevel_2001,giles_multilevel_2008,giles_multilevel_2015} are useful to improve the efficiency of such approaches. While the naive combination of a dynamical low-rank deterministic solver with a Monte Carlo type approach for uncertainty quantification is straightforward, the analysis and optimal design of a more sophisticated scheme requires knowledge of error and cost convergence rates depending on rank. 

Initially, research on using MLMC for stochastic PDE models focused mainly on elliptic problems, notably, Darcy flows with an uncertain diffusion coefficient~\cite{cliffe_multilevel_2011}. However, extensions to an increasing number of other problems exist, such as  parabolic~\cite{barth_multilevel_2013} or hyperbolic problems~\cite{hu_multilevel_2023,mishra_sparse_2012}. 
In particular, MLMC has been applied for the simulation of kinetic equations with particle methods, both at the trajectory level~\cite{lovbak_multilevel_2021,lovbak_accelerated_2023,mortier_multilevel_2022} and the ensemble level~\cite{rosin_multilevel_2014,haji-ali_multilevel_2018,szpruch_iterative_2019}. Methodologically, work on stochastic PDE problems has focused on solving the PDE with different grid resolutions at different levels in the multilevel hierarchy. In this case, it is known that geometric meshes are, in general, close to optimal~\cite{haji-ali_optimization_2016}. However, recent work has advanced MLMC simulations for elliptic problems by integrating hierarchical solvers, an example being multigrid solvers~\cite{kumar_multigrid_2017}, where it has been shown that one can reduce the total computational cost by recycling work from coarser levels~\cite{robbe_recycling_2019}. We also note that multi-index extensions become relevant when considering problems with infinite-dimensional uncertainties, as a hierarchy of approximations is then needed in the random variable~\cite{haji-ali_multi-index_2016}. 

In this work, we combine variance reduction techniques with multi-fidelity dynamical low-rank approximations for the radiation transport equation. We investigate the error of the expectation in terms of both the approximation rank and the spatio-temporal mesh. We then demonstrate how to perform efficient correlated sampling based on the dynamical low-rank method and demonstrate that a control variate can significantly speed up the computation of accurate expectations. The remainder of this paper is structured as follows: In Section~\ref{sec:DLRAMC_DLRA}, we introduce dynamical low-rank approximation \cite{koch_dynamical_2007} as an efficient reduced model-order approach for time-dependent problems and describe the augmented basis update \& Galerkin integrator \cite{ceruti_rank-adaptive_2022}.  In Section~\ref{sec:DLRAMC_UQ}, we extend the well-known robust error bound for the dynamical low-rank approximation \cite{kieri_discretized_2016,ceruti_unconventional_2022} to the probabilistic setting based on the proof from \cite{ceruti_robust_2024} and provide insights into the rank-error relation. We then present a low-rank Monte Carlo and control variate strategy based on these insights. In Section~\ref{sec:DLRAMC_numerical}, we present numerical results for a radiation transport problem that back up our theoretical claims and demonstrates the speedup of our control variate strategy. While our interest lies in higher-dimensional problems, we focus here on a 1D test problem with an uncertain initial condition to limit computational costs. However, our approach straightforwardly applies to other uncertainties and higher-dimensional cases. Finally, in Section~\ref{sec:DLRAMC_conclusion}, we draw conclusions on our results and discuss the potential for a future extension to multilevel Monte Carlo.

%%%%%%%%%%%%%%%%%%%%%%%%%%%%%%%%%%%%%%%%%%%%%%%%%%%%%%%%%%%%%%%%%%%%%%%%%%
\section{Dynamical low-rank approximation} \label{sec:DLRAMC_DLRA}
%%%%%%%%%%%%%%%%%%%%%%%%%%%%%%%%%%%%%%%%%%%%%%%%%%%%%%%%%%%%%%%%%%%%%%%%%%
This section summarizes the fundamentals of dynamical low-rank approximation (DLRA) laid out in \cite{koch_dynamical_2007} and presents the augmented basis-update \& Galerkin \cite{ceruti_rank-adaptive_2022} integrator. We start by assuming a discretization of \eqref{eq:DLRAMC_3Dradiativetransport} in space and direction of travel and write it as
\begin{equation}\label{eq:DLRAMC_MDE}
    \dot{\boldsymbol\Psi}(t) = \textbf{F}(t,\boldsymbol\Psi(t)), \quad
    \boldsymbol{\Psi}(t_0) = \boldsymbol{\Psi_0}\;,
\end{equation}
where $\boldsymbol\Psi\in\R^{m\times n}$, $m,n\gg 1$, is the discretized particle density and $\dot{\boldsymbol\Psi} = \frac{\mathrm{d}}{\mathrm{dt}
}\boldsymbol\Psi$ is the derivative with respect to time. Here, $m$, $n$ represent the number of spatial and directional discretization points, respectively. Then DLRA \cite{koch_dynamical_2007} evolves $\boldsymbol\Psi(t)$ on the low-rank manifold $\mathcal{M}_{r}$, of $m\times n$ rank-$r$ matrices. It does so by projecting the dynamics of the problem onto the tangent space of $\mathcal{M}_{r}$. 

Any given low-rank approximation $\textbf{Y}(t)\in\mathcal{M}_{r}$ to $\boldsymbol\Psi(t)$ can be factorized as
\begin{displaymath}
    \textbf{Y}(t) = \textbf{X}(t)\textbf{S}(t)\textbf{V}(t)^{\top},
\end{displaymath}
where $\textbf{X}(t)\in\R^{m\times r}$, $\textbf{V}(t)\in\R^{n\times r}$ have orthonormal columns and $\textbf{S}(t)\in\R^{r\times r}$ is invertible. Let $\mathcal{T}_{\textbf{Y}}\mathcal{M}_{r}$ denote the tangent space of the low-rank manifold $\mathcal{M}_{r}$ at $\textbf{Y}$. To ensure that the low-rank approximation $\textbf{Y}(t)$ remains of rank $r$ at any given time, DLRA solves the optimization problem
\begin{displaymath}
    \underset{\dot{\textbf{Y}}(t)\in\mathcal{T}_{\textbf{Y}(t)}\mathcal{M}_{r}}{\min}\lVert \dot{\textbf{Y}}(t) - \textbf{F}(t,\textbf{Y}(t)) \rVert,
\end{displaymath}
where $\lVert \cdot \rVert$ is the Frobenius norm. This distance is minimized by the orthogonal projection of the right-hand side of \eqref{eq:DLRAMC_MDE} onto the tangent space $\mathcal{T}_{\textbf{Y}(t)}\mathcal{M}_{r} $, denoted by $\textbf{P}(\textbf{Y}(t))$, \cite[Lemma~4.1]{koch_dynamical_2007}, i.e. 
\begin{displaymath}
    \dot{\textbf{Y}} = \boldsymbol{\mathcal{P}}(\textbf{Y}(t))\textbf{F}(t,\textbf{Y}(t)),
\qquad
\text{where}
\qquad
    \boldsymbol{\mathcal{P}}(\textbf{Y})\textbf{Z} = \textbf{X}\textbf{X}^{\top}\textbf{Z} - \textbf{X}\textbf{X}^{\top}\textbf{Z}\textbf{V}\textbf{V}^{\top} +\textbf{Z}\textbf{V}\textbf{V}^{\top}.
\end{displaymath}
Thus the evolution equations for the factors ensuring $\textbf{Y}(t)\in\mathcal{M}_{r}$ for all $t$, are \cite{koch_dynamical_2007}
\begin{align*}
    \dot{\textbf{S}}(t) &= \textbf{X}(t)^{\top}\textbf{F}(t,\textbf{Y}(t))\textbf{V}(t),\\
    \dot{\textbf{X}}(t) &= (\textbf{I}-\textbf{X}(t)\textbf{X}(t)^{\top})\textbf{F}(t,\textbf{Y}(t))\textbf{V}(t)\textbf{S}(t)^{-1},\\
    \dot{\textbf{V}}(t) &= (\textbf{I} - \textbf{V}(t)\textbf{V}(t)^{\top})\textbf{F}(t,\textbf{Y}(t))^{\top}\textbf{X}(t)\textbf{S}(t)^{-\top}.
\end{align*}
Since the rank of the solution is not known beforehand, it is often over-approximated in practice. This introduces several near-zero singular values into the approximation, causing conventional time integrators to fail in achieving convergence. Several robust approaches have been developed to counter this problem, such as the projector-splitting \cite{lubich_projector-splitting_2014} or basis-update \& Galerkin integrators \cite{ceruti_unconventional_2022,ceruti_rank-adaptive_2022,ceruti_parallel_2024} and projection-based methods \cite{kieri_projection_2019}. In this work, we use the augmented BUG integrator \cite{ceruti_rank-adaptive_2022} due to its favorable stability properties for radiation transport \cite{kusch_stability_2023}.

Let $\textbf{Y}_{k}:=\textbf{Y}(t_{k}),\textbf{X}_{k}:=\textbf{X}(t_{k}),\textbf{S}_{k}:=\textbf{S}(t_{k}),\textbf{V}_{k}:=\textbf{V}(t_{k})$, where $t_{k} = t_{0} + k h$ for some starting time $t_{0}\in\R_{\geq0}$ and step size $h>0$. Then the rank-$r$ approximation at $t_{0}$ is $\textbf{Y}_{0} = \textbf{X}_{0}\textbf{S}_{0}\textbf{V}_{0}^{\top} $. In the following, we outline one update step according to \cite{ceruti_rank-adaptive_2022}, where the factors $\textbf{X}_{0}$, $\textbf{S}_{0}$, $\textbf{V}_{0}$ are updated to $\textbf{X}_{1}$, $\textbf{S}_{1}$, $\textbf{V}_{1}$ at time $t_1$ in three sub-steps. Since the goal is to gain a better understanding of uncertainty propagation in the low-rank framework, the rank of the approximation is not adapted according to a truncation tolerance, as proposed in \cite{ceruti_rank-adaptive_2022}, but truncated to a fixed rank $r$. 
\begin{enumerate}
    \item \textbf{Basis augmentation}\\
        \textbf{K-step}: For $\textbf{K}(t) = \textbf{X}(t)\textbf{S}(t)$, integrate from $t_{0}$ to $t_{1}$
        \begin{displaymath}
        \dot{\textbf{K}}(t) = \textbf{F}(t,\textbf{K}(t)\textbf{V}_{0}^{\top})\textbf{V}_{0}, \quad \textbf{K}(t_{0}) = \textbf{X}_{0}\textbf{S}_{0},
        \end{displaymath}
        compute $\widehat{\textbf{X}}\in\R^{m\times 2r}$ as an orthonormal basis of $[\textbf{K}(t_{1}), \textbf{X}_{0}]$ and store $\textbf{M} = \widehat{\textbf{X}}^{\top}\textbf{X}_{0} $.
        \item[] \textbf{L-step}: For $\textbf{L}(t) = \textbf{V}(t)\textbf{S}(t)^{\top} $, integrate from $t_{0}$ to $t_{1}$
        \begin{displaymath}
        \dot{\textbf{L}}(t) = \textbf{F}(t,\textbf{X}_{0}\textbf{L}(t)^{\top})^{\top}\textbf{X}_{0}, \quad \textbf{L}(t_{0}) = \textbf{V}_{0}\textbf{S}_{0}^{\top},
        \end{displaymath}
        compute $\widehat{\textbf{V}}\in\R^{n\times 2r}$ as an orthonormal basis of $[\textbf{L}(t_{1}), \textbf{V}_{0}]$ and store $\textbf{N} = \widehat{\textbf{V}}^{\top}\textbf{V}_{0} $.
    \item \textbf{Galerkin step}\\
        \textbf{S-step}: Integrate from $t_{0}$ to $t_{1}$
        \begin{equation}\label{eq:DLRAMC_Sstep}
        \dot{\widehat{\textbf{S}}}(t) = \widehat{\textbf{X}}^{\top}\textbf{F}(t,\widehat{\textbf{X}}\widehat{\textbf{S}}(t)\widehat{\textbf{V}}^{\top})\widehat{\textbf{V}}, \quad \widehat{\textbf{S}}(t_{0}) = \textbf{M}\textbf{S}_{0}\textbf{N}^{\top}.
        \end{equation}
    \item \textbf{Truncation to rank $r$}\\
    Compute the singular value decomposition of $\widehat{\textbf{S}}(t_{1}) = \textbf{P}\boldsymbol{\Sigma}\textbf{Q}^{\top}$. Let $\mathbf{S}_{1}\in\mathbb{R}^{r\times r}$ be the diagonal matrix with the $r$ largest singular values of $\boldsymbol{\Sigma}$ and set $\mathbf{P}_1\in\mathbb{R}^{2r\times r}$ and $\mathbf{Q}_{1}\in\mathbb{R}^{2r\times r}$ contain the first $r$ columns of $\textbf{P}$ and $\textbf{Q}$.
    Then, $\textbf{X}_{1} = \widehat{\textbf{X}}\mathbf{P}_{1}$ and $\textbf{V}_{1} = \widehat{\textbf{V}}\mathbf{Q}_{1}$.
\end{enumerate}
Finally, the approximation at $t_{1}$ is set as $\textbf{Y}_{1} = \textbf{X}_{1}\textbf{S}_{1}\textbf{V}_{1}^{\top} $. Note, if $ \widehat{\textbf{Y}}(t) = \widehat{\textbf{X}}\widehat{\textbf{S}}(t)\widehat{\textbf{V}}^{\top} $ the augmented solution is obtained by solving the following reformulation of \eqref{eq:DLRAMC_Sstep}:
\begin{equation}\label{eq:DLRAMC_SstepRe}
    \dot{\widehat{\textbf{Y}}}(t) = \widehat{\textbf{X}}\widehat{\textbf{X}}^{\top}\textbf{F}(t,\widehat{\textbf{Y}}(t))\widehat{\textbf{V}}\widehat{\textbf{V}}^{\top}, \quad \widehat{\textbf{Y}}(t_{0}) = \textbf{Y}_{0}.
\end{equation}

%%%%%%%%%%%%%%%%%%%%%%%%%%%%%%%%%%%%%%%%%%%%%%%%%%%%%%%%%%%%%%%%%%%%%%%%%%
% \subsection{BUG integrator}
%%%%%%%%%%%%%%%%%%%%%%%%%%%%%%%%%%%%%%%%%%%%%%%%%%%%%%%%%%%%%%%%%%%%%%%%%%

%%%%%%%%%%%%%%%%%%%%%%%%%%%%%%%%%%%%%%%%%%%%%%%%%%%%%%%%%%%%%%%%%%%%%%%%%
% \subsection{Augmented BUG}
%%%%%%%%%%%%%%%%%%%%%%%%%%%%%%%%%%%%%%%%%%%%%%%%%%%%%%%%%%%%%%%%%%%%%%%%%%

%%%%%%%%%%%%%%%%%%%%%%%%%%%%%%%%%%%%%%%%%%%%%%%%%%%%%%%%%%%%%%%%%%%%%%%%%%
\section{Low-rank estimators} \label{sec:DLRAMC_UQ}
%%%%%%%%%%%%%%%%%%%%%%%%%%%%%%%%%%%%%%%%%%%%%%%%%%%%%%%%%%%%%%%%%%%%%%%%%%
In this section, we introduce uncertainty to \eqref{eq:DLRAMC_MDE} and investigate the propagation of uncertainty through the solution. We do so by constructing Monte Carlo and control variate estimators based on low-rank approximation of the solution. Furthermore, we show that when used together with Monte Carlo sampling, the augmented BUG integrator is robust to the presence of small singular values.

Let $\nu$ be a scalar random variable with probability density function $p(\nu)$. Then if $\boldsymbol{\Psi}(t;\nu) \in \R^{m\times n} $ is the discretized particle density subject to the random variable $\nu$, the uncertain matrix differential equation (MDE) reads
\begin{equation}\label{eq:DLRAMC_UncertainMDE}
    \dot{\boldsymbol{\Psi}}(t;\nu) = \textbf{F}(t,\boldsymbol{\Psi}(t;\nu)), \qquad \boldsymbol{\Psi}(t_{0};\nu) = \boldsymbol{\Psi}_{0}(\nu).
\end{equation}
For simplicity, we choose the uncertain initial conditions that, for example, arise from uncertainty in the position or strength of particle sources in \eqref{eq:DLRAMC_3Dradiativetransport}~\cite{kusch_low-rank_2022}. However, the methods developed in this section can be analogously applied to other types of uncertainties. We assume that, for a fixed $\nu$, \eqref{eq:DLRAMC_UncertainMDE} has a low-rank solution which is a reasonable assumption, for instance for the radiative transfer equation~\cite{einkemmer_low-rank_2018,kusch_low-rank_2022,kusch_robust_2023}.
% For example, if the uncertainty is in the position or strength of the particle source, it is likely that this will not severely affect the rank of the solution. 

We are interested in a lower-dimensional function of the particle density, such as the scalar flux, at a given time $t_{k} >0$. We therefore define the map $\nu\to \mathcal{G}(\boldsymbol{\Psi}(t = t_{k};\nu)) \in\R^{\widetilde{m}\times \widetilde{n}} $, where $\widetilde{m},\widetilde{n}\in\N$. Then the quantity of interest (QoI), $Q\in\R^{\widetilde{m}\times \widetilde{n}}$, is the expectation of $\mathcal{G}(\boldsymbol{\Psi}(t_{k};\nu))$. The expectation and variance of $\mathcal{G}(\boldsymbol{\Psi}(t_{k};\nu))$ are, respectively,
\begin{displaymath}
    Q \coloneqq\mathbb{E}_{p}\left[\mathcal{G}(\boldsymbol{\Psi}(t_{k};\nu))\right] = \int_{-\infty}^\infty\mathcal{G}(\boldsymbol{\Psi}(t_{k};\nu)) p(\nu)~\mathrm{d}\nu \qquad \text{and}
\end{displaymath}
\begin{displaymath}
    \mathrm{Var}(\mathcal{G}(\boldsymbol{\Psi}(t_{k};\nu))) = \int_{-\infty}^\infty\left\lVert \mathcal{G}(\boldsymbol{\Psi}(t_{k};\nu)) - \mathbb{E}_{p}\left[\mathcal{G}(\boldsymbol{\Psi}(t_{k};\nu))\right] \right\rVert^{2} p(\nu)~\mathrm{d}\nu,
\end{displaymath}
where $\lVert \cdot \rVert$ is a matrix norm induced by a suitable inner product $\langle \cdot,\cdot\rangle$, e.g., the Frobenius norm induced by the Frobenius inner product.

In Section~\ref{sec:DLRAMC_MonteCarloEst} we present a low-rank Monte Carlo estimator to $Q$ and provide error bounds on the low-rank approximation under Monte Carlo sampling. In Section~\ref{sec:DLRAMC_CVEst}, we introduce a control variate strategy for variance reduction.

%%%%%%%%%%%%%%%%%%%%%%%%%%%%%%%%%%%%%%%%%%%%%%%%%%%%%%%%%%%%%%%%%%%%%%%%%%
\subsection{Low-rank Monte Carlo estimator} \label{sec:DLRAMC_MonteCarloEst}
%%%%%%%%%%%%%%%%%%%%%%%%%%%%%%%%%%%%%%%%%%%%%%%%%%%%%%%%%%%%%%%%%%%%%%%%%%
To estimate $Q$, we need to know the exact solution, $\boldsymbol{\Psi}(t;\nu)$, to \eqref{eq:DLRAMC_UncertainMDE} at $t = t_{k}$. When the exact solution is unknown it is approximated with a standard time-stepping scheme like the explicit Euler or a higher-order Runge-Kutta method. Let $\widetilde{\boldsymbol{\Psi}}_{k}(\nu) \coloneqq \widetilde{\boldsymbol{\Psi}}(t_{k};\nu)$ denote the approximation to $\boldsymbol{\Psi}_{k}(\nu)\coloneqq\boldsymbol{\Psi}(t=t_{k};\nu)$ computed with a time-stepping scheme. Since $\widetilde{\boldsymbol{\Psi}}_{k}(\nu) \approx \boldsymbol{\Psi}_{k}(\nu)$, we approximate $Q$ by $ \mathbb{E}_{p}[\mathcal{G}(\widetilde{\boldsymbol{\Psi}}_{k}(\nu))] \eqqcolon Q_{\mathrm{full}}$ and use the Monte Carlo method to construct the estimator, $\widehat{Q}_{\mathrm{full};\mathrm{MC}}$. However, this estimator is computationally expensive and memory-intensive since it requires solving and storing the solution of a high-dimensional MDE for each realization of the random variable. 

To overcome the computational inefficiency of $\widehat{Q}_{\mathrm{full};\mathrm{MC}}$, we make use of the assumption that \eqref{eq:DLRAMC_UncertainMDE} has a low-rank solution. As shown in Section~\ref{sec:DLRAMC_DLRA}, when the underlying solution to a time-dependent problem has a low-rank structure, DLRA~\cite{koch_dynamical_2007} can be used to approximate the solution. Let $\textbf{Y}_{k}(\nu)=\textbf{X}_{k}(\nu)\textbf{S}_{k}(\nu)\textbf{V}_{k}(\nu)^{\top} $ denote the rank-$r$ approximation to $\boldsymbol{\Psi}_{k}(\nu)$ computed using the augmented BUG integrator described in Section~\ref{sec:DLRAMC_DLRA}. Then instead of approximating $Q$ we estimate $Q_{r} \coloneqq \mathbb{E}_{p}[\mathcal{G}(\textbf{Y}_{k}(\nu))]$ which is computationally efficient and has lower memory requirements.

Let $\nu_1,...,\nu_N \sim p(\nu)$ be $N$ independent and identically distributed random realizations of $\nu$. Then we define the low-rank Monte Carlo estimator as
\begin{displaymath}
 \widehat{Q}_{r;\mathrm{MC}} := \frac{1}{N} \sum_{i=1}^{N} \mathcal{G}(\textbf{Y}_{k}(\nu_{i})) .
\end{displaymath}
The mean squared error (MSE) of this estimator is
\begin{align}
    \mathrm{MSE}(\widehat{Q}_{r;\mathrm{MC}}) &= \mathbb{E}_{p}\left[\left\lVert \frac{1}{N} \sum_{i=1}^{N} \mathcal{G}(\textbf{Y}_{k}(\nu_{i})) - \mathbb{E}_{p}[\mathcal{G}(\boldsymbol{\Psi}_{k}(\nu))] \right\rVert^{2}\right]\notag\\
    &= \mathbb{E}_{p}\left[\left\lVert \widehat{Q}_{r;\mathrm{MC}} - \mathbb{E}_{p}[\widehat{Q}_{r;\mathrm{MC}}]\right\rVert^{2}\right] + \left\lVert \mathbb{E}_{p}[\widehat{Q}_{r;\mathrm{MC}}] - Q\right\rVert^{2}.\label{eq:DLRAMC_MSE}
\end{align}
Thus, the MSE of the estimator separates into the variance of the Monte Carlo estimator and a bias term, describing errors caused by the model, numerical discretization, and low-rank approximation. The variance can be written as
\begin{displaymath}\label{eq:DLRAMC_MCvariance}
    \mathbb{E}_{p}\Big[\left\lVert \widehat{Q}_{r;\mathrm{MC}} - \mathbb{E}_{p}[\widehat{Q}_{r;\mathrm{MC}}]\right\rVert^{2}\Big] =  \mathrm{Var}\Big[ \widehat{Q}_{r;\mathrm{MC}} \Big] = \frac{1}{N} \mathrm{Var}\Big[ \mathcal{G}(\textbf{Y}_{k}(\nu))\Big].
\end{displaymath}
Similarly, the bias can be written as
\begin{displaymath}
    \left\lVert\mathbb{E}_{p}[\widehat{Q}_{r;\mathrm{MC}}] - Q\right\rVert = \left\lVert Q_{r} - Q\right\rVert = \left\lVert \mathbb{E}_{p}[\mathcal{G}(\textbf{Y}_{k}(\nu)) - \mathcal{G}(\boldsymbol{\Psi}_{k}(\nu))]\right\rVert,
\end{displaymath}
where we use $ \mathbb{E}_{p}[\widehat{Q}_{r;\mathrm{MC}}] = Q_{r} $. The bias further splits into contributions from the low-rank approximation and those from the underlying model and its discretization. We now provide a robust error bound for the low-rank approximation in the presence of small singular values due to over-approximation of rank, i.e. for the error contribution due to low-rank approximation under uncertainty.  

In the deterministic setting, the error bound of \cite{ceruti_rank-adaptive_2022} shows that the augmented BUG integrator is robust to small singular values. However, this result does not hold in the presence of uncertainty. The next theorem shows that the augmented BUG integrator applied to the uncertain matrix-valued differential equation \eqref{eq:DLRAMC_MDE} is robust to over-approximation and how the bounding constants depend on the random variable. The proof closely follows the local error bound for the mid-point BUG integrator~\cite[Theorem~1]{ceruti_robust_2024} and uses the following assumptions, where $\lVert \cdot \rVert = \lVert \cdot \rVert_{F}$ is the Frobenius norm:
 \begin{enumerate}
        \item[A.1] $\normalfont\textbf{F}$, the right-hand side of \eqref{eq:DLRAMC_MDE}, is Lipschitz-continuous, bounded, and independent of $\nu$: for all $\normalfont\textbf{Y},\widetilde{\textbf{Y}}\in\mathbb{R}^{m\times n}$ and $t_{0} \leq t \leq t_{k}$, 
        \begin{displaymath}
             \normalfont \lVert\textbf{F}(t,\textbf{Y}) - \textbf{F}(t,\widetilde{\textbf{Y}}) \rVert \leq L\lVert \textbf{Y} - \widetilde{\textbf{Y}} \rVert, \qquad \lVert \textbf{F}(t,\textbf{Y}) \rVert \leq B.
        \end{displaymath}

        \item[A.2] For $ \normalfont \textbf{F}(t,\textbf{Y}(t;\nu)) = \textbf{M}(t,\textbf{Y}(t;\nu)) + \textbf{R}(t,\textbf{Y}(t;\nu)) $, where $\normalfont\textbf{M}(t,\textbf{Y}(t;\nu)) $ \\$= \boldsymbol{\mathcal{P}}(\textbf{Y}(t;\nu))\textbf{F}(t,\textbf{Y}(t;\nu))$. The non-tangential part $\normalfont \textbf{R}(t,\textbf{Y}(t;\nu))$ is $\varepsilon$-small:
    \begin{displaymath}
           \normalfont  \lVert \textbf{R}(t,\textbf{Y}(t;\nu)) \rVert \leq c_{1}(\nu)\varepsilon, \quad c_{1} \in L^1.
        \end{displaymath}
        \item[A.3] The error in the initial value is $\delta$-small:
        \begin{displaymath}
             \normalfont\lVert \textbf{Y}_{0}(\nu) - \boldsymbol{\Psi}_{0}(\nu) \rVert \leq c_{0}(\nu)\delta, \quad c_{0} \in L^1.
        \end{displaymath}
    \end{enumerate}

\begin{theorem}[Local error bound]\label{theorem:DLRAMC_DetRobErrBd}
    Let $\normalfont\boldsymbol{\Psi}(t;\nu)$ denote the solution of the uncertain matrix differential equation \eqref{eq:DLRAMC_MDE} and $ \normalfont\textbf{Y}_{1}(\nu)$ denote the rank-$r$ approximation to $ \boldsymbol{\Psi}_{1}(\nu)\coloneqq \boldsymbol{\Psi}(t_{1};\nu) $ at $t_{1} = t_{0} + h$ obtained after one step of the augmented BUG integrator with step-size $h>0$ and truncation error $\vartheta_{1}(\nu)$. Then, if the assumptions (A.1-A.2) are satisfied, we have the following local error bound on the expected value of the solution
    \begin{displaymath}
       \normalfont \left\lVert \mathbb{E}_{p}[ \textbf{Y}_{1}(\nu) - \boldsymbol{\Psi}_{1}(\nu) ] \right\rVert \leq \mathbb{E}_{p}[\widetilde{c}_{1}(\nu)]\varepsilon h + 10LBh^{2}\,,
    \end{displaymath}
    where $\widetilde{c}_{1}$ depends on $c_{1}$, $L$, and $t_{\mathrm{end}}$.
\end{theorem}
\begin{remark}
    We fix some notation used in the proof. Let $\textbf{Y}_{0}^{\nu}=\textbf{X}_{0}^{\nu}\textbf{S}_{0}^{\nu}\textbf{V}_{0}^{\nu,\top}$ be a realization of $\textbf{Y}_{0}(\nu)$, i.e. $\textbf{Y}_{0}^{\nu}\sim\textbf{Y}_{0}(\nu)$. Then the exact solution at time $t_{1} = t_{0} + h$, with the initial condition $\textbf{Y}_{0}^{\nu}$, is denoted by $\boldsymbol{\Psi}_{1}^{\nu}\sim \boldsymbol{\Psi}_{1}(\nu)$. Similarly, $\textbf{Y}_{1}^{\nu}=\textbf{X}_{1}^{\nu}\textbf{S}_{1}^{\nu}\textbf{V}_{1}^{\nu,\top}\sim \textbf{Y}_{1}({\nu})$ denotes the low-rank approximation at $t_{1}$ computed with the augmented BUG integrator and $\vartheta_1^{\nu}\sim \vartheta_1(\nu)$ denotes the truncation error made at $t_1$. Additionally, let $\widehat{\textbf{Y}}_{1}^{\nu}=\widehat{\textbf{X}}^{\nu}\widehat{\textbf{S}}^{\nu}\widehat{\textbf{V}}^{\nu,\top}\sim \widehat{\textbf{Y}}_{1}^{\nu}$ denote the augmented solution before truncation.
\end{remark}
\begin{proof}
% Let $\textbf{X}_{i}^{\nu}:=\textbf{X}(t_i;\nu),\textbf{V}_{i}^{\nu}:=\textbf{V}(t_i;\nu),\; i\in\{0,1\}$ and $\widehat{\textbf{X}}^{\nu} \in \mathbb{R}^{m\times 2r},\widehat{\textbf{V}}^{\nu} \in \mathbb{R}^{n\times 2r}$ denote the chosen augmented bases. \\
    We start by using Jensen's and Cauchy-Schwarz (CS) inequality to get,
    \begin{displaymath}
        \left\lVert \mathbb{E}_{p}[ \textbf{Y}_{1}(\nu) - \boldsymbol{\Psi}_{1}(\nu)  ] \right\rVert \leq \mathbb{E}_{p}\left[ \lVert \widehat{\textbf{Y}}_{1}(\nu) - \boldsymbol{\Psi}_{1}(\nu) \rVert \right] + \mathbb{E}_{p}\left[ \lVert \widehat{\textbf{Y}}_{1}(\nu) - \textbf{Y}_{1}(\nu) \rVert \right].
    \end{displaymath}
    % where $\lVert \boldsymbol{\Psi}_{1}^{\nu} - \textbf{Y}_{1}^{\nu} \rVert$ is a realization of $\lVert \boldsymbol{\Psi}_{1}(\nu) - \textbf{Y}_{1}(\nu) \rVert$ \pia{$\lVert \boldsymbol{\Psi}_{1}(\nu) - \textbf{Y}_{1}(\nu) ] \rVert$? and I'm not sure about the phrasing with "expectation is over all such realizations" to me it is already a bit like talking about a MC approx. maybe its better to refer to the pdf instead, but I'm also not so confident on this, maybe @Sebastian or @Emil can look as well} and the expectation is over all such realizations. 
    For each realization we have
    \begin{displaymath}
        \lVert \widehat{\textbf{Y}}_{1}^{\nu} - \boldsymbol{\Psi}_{1}^{\nu} \rVert \overset{\small{\mathrm{CS}}}{\leq} \lVert \boldsymbol{\Psi}_{1}^{\nu} - \widehat{\textbf{X}}^{\nu}\widehat{\textbf{X}}_{1}^{\nu,\top}\boldsymbol{\Psi}_{1}^{\nu}\widehat{\textbf{V}}^{\nu}\widehat{\textbf{V}}^{\nu,\top} \rVert + \lVert \widehat{\textbf{X}}^{\nu}\widehat{\textbf{X}}^{\nu,\top}\boldsymbol{\Psi}_{1}^{\nu}\widehat{\textbf{V}}^{\nu}\widehat{\textbf{V}}^{\nu,\top} - \widehat{\textbf{Y}}_{1}^{\nu} \rVert.
    \end{displaymath}
    Adding and subtracting $\widehat{\textbf{X}}^{\nu}\widehat{\textbf{X}}^{\nu,\top}\boldsymbol{\Psi}_{1}^{\nu}$ from the first term on the right-hand side yields
    \begin{equation}\label{eq:DLRAMC_LEeq2}
        \Vert \boldsymbol{\Psi}_{1}^{\nu} - \widehat{\textbf{X}}^{\nu}\widehat{\textbf{X}}^{\nu,\top}\boldsymbol{\Psi}_{1}^{\nu}\widehat{\textbf{V}}^{\nu}\widehat{\textbf{V}}^{\nu,\top} \rVert \leq \lVert \boldsymbol{\Psi}_{1}^{\nu} - \widehat{\textbf{X}}^{\nu}\widehat{\textbf{X}}^{\nu,\top}\boldsymbol{\Psi}_{1}^{\nu} \rVert + \lVert \boldsymbol{\Psi}_{1}^{\nu} - \boldsymbol{\Psi}_{1}^{\nu}\widehat{\textbf{V}}^{\nu}\widehat{\textbf{V}}^{\nu,\top} \rVert.
    \end{equation}
    From the definition of $\widehat{\textbf{X}}^{\nu}$, $ (\textbf{I} - \widehat{\textbf{X}}^{\nu}\widehat{\textbf{X}}^{\nu,\top})\textbf{K}^{\nu}(t_{1})\textbf{V}_{0}^{\nu,\top} = 0 $, we write the first term of the right-hand side of \eqref{eq:DLRAMC_LEeq2} as
    \begin{align*}
        \lVert \boldsymbol{\Psi}_{1}^{\nu} - \widehat{\textbf{X}}^{\nu}\widehat{\textbf{X}}^{\nu,\top}\boldsymbol{\Psi}_{1}^{\nu} \rVert &= \lVert (\textbf{I} - \widehat{\textbf{X}}^{\nu}\widehat{\textbf{X}}^{\nu,\top})(\boldsymbol{\Psi}_{1}^{\nu} - \textbf{K}^{\nu}(t_{1})\textbf{V}_{0}^{\nu,\top}) \rVert\\
        &\!\!\!\!\!\!\overset{\small{\mathrm{FTC}}}{\leq} \int_{t_{0}}^{t_{1}}\lVert (\textbf{I} - \widehat{\textbf{X}}^{\nu}\widehat{\textbf{X}}^{\nu,\top})( \textbf{F}(t,\boldsymbol{\Psi}^{\nu}(t)) - \textbf{F}(t,\textbf{K}^{\nu}(t)\textbf{V}_{0}^{\nu,\top})\textbf{V}_{0}^{\nu}\textbf{V}_{0}^{\nu,\top}) \rVert\hspace{1mm}\mathrm{dt},
    \end{align*}
    where FTC stands for the fundamental theorem of calculus. We can write
    \begin{align*}
        \textbf{F}(t,\boldsymbol{\Psi}^{\nu}(t)) &=  \textbf{F}(t_{0},\textbf{Y}^{\nu}_{0}) + \textbf{F}(t,\boldsymbol{\Psi}^{\nu}(t)) - \textbf{F}(t_{0},\textbf{Y}^{\nu}_{0})\\
        &= \textbf{M}(t_{0},\textbf{Y}^{\nu}_{0}) +  \textbf{F}(t,\boldsymbol{\Psi}^{\nu}(t)) - \textbf{F}(t_{0},\textbf{Y}^{\nu}_{0}) + \textbf{R}(t_{0},\textbf{Y}^{\nu}_{0}),
    \end{align*}
    \begin{align*}
        \text{and} \quad \textbf{F}(t,\textbf{K}^{\nu}(t)\textbf{V}_{0}^{\nu,\top})\textbf{V}_{0}^{\nu}\textbf{V}_{0}^{\nu,\top} &= \textbf{F}(t_{0},\textbf{Y}_{0}^{\nu})\textbf{V}_{0}^{\nu}\textbf{V}_{0}^{\nu,\top} \\ & \qquad  + \left( \textbf{F}(t,\textbf{K}^{\nu}(t)\textbf{V}_{0}^{\nu,\top})\textbf{V}_{0}^{\nu}\textbf{V}_{0}^{\nu,\top} - \textbf{F}(t_{0},\textbf{Y}_{0}^{\nu})\textbf{V}_{0}^{\nu}\textbf{V}_{0}^{\nu,\top} \right).
    \end{align*}
    Since $\lVert \textbf{I} - \widehat{\textbf{X}}^{\nu}\widehat{\textbf{X}}^{\nu,\top} \rVert\leq 1$, we have 
    % \jonas{Where did the $M$ and $FV_0V_0^T$ go? This does not work.}\chins{@Jonas is this comment still valid with the augmented BUG?}
    \begin{align*}
            \lVert \boldsymbol{\Psi}_{1}^{\nu} - \widehat{\textbf{X}}^{\nu}\widehat{\textbf{X}}^{\nu,\top}\boldsymbol{\Psi}_{1}^{\nu} \rVert &\overset{\mathrm{CS}}{\leq} \int_{t_{0}}^{t_{1}}\lVert \textbf{F}(t,\textbf{K}^{\nu}(t)\textbf{V}_{0}^{\nu,\top}) - \textbf{F}(t_{0},\textbf{Y}_{0}^{\nu}) \rVert \hspace{1mm}\mathrm{dt} \\ &\qquad + \int_{t_{0}}^{t_{1}}\lVert \textbf{F}(t,\boldsymbol{\Psi}^{\nu}(t)) - \textbf{F}(t_{0},\textbf{Y}_{0}^{\nu}) \rVert \hspace{1mm}\mathrm{dt} + \int_{t_{0}}^{t_{1}}\lVert \textbf{R}(t_{0},\textbf{Y}^{\nu}_{0}) \rVert \hspace{1mm}\mathrm{dt}\\
            &\!\!\!\!\!\!\!\!\!\!\!\overset{A.1, A.2}{\leq} L\int_{t_{0}}^{t_{1}}\lVert \textbf{K}^{\nu}(t)\textbf{V}_{0}^{\nu,\top} - \textbf{Y}_{0}^{\nu} \rVert \hspace{1mm}\mathrm{dt} + L\int_{t_{0}}^{t_{1}}\lVert \boldsymbol{\Psi}^{\nu}(t) - \textbf{Y}_{0}^{\nu} \rVert \hspace{1mm}\mathrm{dt} + c_{1}^{\nu}h\varepsilon\,.
    \end{align*}
   Thus using FTC and $A.1$ we get,
    \begin{displaymath}
        \lVert \boldsymbol{\Psi}_{1}^{\nu} - \widehat{\textbf{X}}^{\nu}\widehat{\textbf{X}}^{\nu,\top}\boldsymbol{\Psi}_{1}^{\nu} \rVert \leq 2LBh^{2} + c_{1}^{\nu}h\varepsilon.
    \end{displaymath}
    Using $ (\textbf{I} - \widehat{\textbf{V}}^{\nu}\widehat{\textbf{V}}^{\nu,\top})\textbf{L}^{\nu}(t_{1})\textbf{X}_{0}^{\nu,\top} = 0 $ we get an analogous expression for $\lVert \boldsymbol{\Psi}_{1}^{\nu} - \boldsymbol{\Psi}_{1}^{\nu}\widehat{\textbf{V}}^{\nu}\widehat{\textbf{V}}^{\nu,\top} \rVert$. Thus, \eqref{eq:DLRAMC_LEeq2} becomes
    \begin{displaymath}
        \Vert \boldsymbol{\Psi}_{1}^{\nu} - \widehat{\textbf{X}}^{\nu}\widehat{\textbf{X}}^{\nu,\top}\boldsymbol{\Psi}_{1}^{\nu}\widehat{\textbf{V}}^{\nu}\widehat{\textbf{V}}^{\nu,\top} \rVert \leq 4LBh^{2} + 2c_{1}^{\nu}h\varepsilon.
    \end{displaymath}
    Now since $\widehat{\textbf{Y}}^{\nu}(t) = \widehat{\textbf{X}}^{\nu}\widehat{\textbf{S}}^{\nu}(t)\widehat{\textbf{V}}^{\nu}  $ we have $\Vert  \widehat{\textbf{X}}^{\nu} \widehat{\textbf{X}}^{\nu,\top}\boldsymbol{\Psi}_1^{\nu}\widehat{\textbf{V}}^{\nu} \widehat{\textbf{V}}^{\nu,\top} - \widehat{\textbf{Y}}_{1}^{\nu} \Vert$
    \begin{align*}
        &\overset{\mathrm{FTC}}{\leq}\, \int_{t_{0}}^{t_{1}}\Vert \mathbf{F}(t,\boldsymbol{\Psi}^{\nu}(t)) - \mathbf{F}(t,\widehat{\textbf{Y}}^{\nu}(t)) \Vert\,\hspace{1mm}\mathrm{dt} \\
        &\overset{A.1}{\leq}\, L\int_{t_{0}}^{t_{1}}\Vert \boldsymbol{\Psi}^{\nu}(t) - \widehat{\textbf{Y}}^{\nu}(t) \Vert~\mathrm{dt}\\
        &\overset{\mathrm{FTC}}{\leq}\, L\int_{t_{0}}^{t_{1}}\int_{t_{0}}^{t}\Vert \mathbf{F}(s,\boldsymbol{\Psi}^{\nu}(s)) - \widehat{\textbf{X}}^{\nu}\widehat{\textbf{X}}^{\nu,\top}\mathbf{F}(s,\widehat{\textbf{Y}}^{\nu}(s))\widehat{\textbf{V}}^{\nu}\widehat{\textbf{V}}^{\nu,\top} \Vert\,\hspace{1mm}\mathrm{ds}\,\mathrm{dt}\\
        &\!\!\!\overset{\mathrm{CS},A.1}{\leq}\, LBh^2\,.
    \end{align*}
    Putting it all together yields
    \begin{displaymath}
        \Vert \widehat{\textbf{Y}}_1^{\nu} - \boldsymbol{\Psi}_1^{\nu} \Vert \leq 5LBh^2 +   2c_{1}^{\nu}\varepsilon h \,.
    \end{displaymath}
     The error made by truncating to rank-$r$ is given by \cite[Lemma~3]{kieri_projection_2019}: In the Frobenius norm we have
    \begin{displaymath}
        \begin{aligned}
            \lVert \widehat{\textbf{Y}}_1^{\nu} - \textbf{Y}_1^{\nu} \rVert &= \underset{\textbf{Z}\in\mathcal{M}_{r}}{\mathrm{min}} \lVert \widehat{\textbf{Y}}_1^{\nu} - \textbf{Z} \rVert\\
            &\leq   \lVert \widehat{\textbf{Y}}_1^{\nu} - \boldsymbol{\Psi}_1^{\nu} \rVert + \underset{\textbf{Z}\in\mathcal{M}_{r}}{\mathrm{min}}\lVert \boldsymbol{\Psi}_1^{\nu} -  \textbf{Z} \rVert + \leq 5LBh^2 +C_{1}^{\nu}\varepsilon h \,,
        \end{aligned}
    \end{displaymath}
    where the bound for the second term, $\underset{\textbf{Z}\in\mathcal{M}_{r}}{\mathrm{min}}\lVert \boldsymbol{\Psi}_1^{\nu} -  \textbf{Z} \rVert = \mathcal{O}(\varepsilon h)$, is obtained from \cite[Lemma~1]{kieri_projection_2019} and $C_{1}^{\nu}$ depends on $L$ and $t_{\mathrm{end}}$. Thus, we get the following local error bound
    \begin{displaymath}
        \left\lVert \mathbb{E}_{p}[ \textbf{Y}_{1}(\nu) - \boldsymbol{\Psi}_{1}(\nu) ] \right\rVert \leq  \mathbb{E}_{p}[\widetilde{c}_{1}(\nu)]\varepsilon h + 10LBh^{2}.
    \end{displaymath}

\end{proof}

\begin{theorem}[Global error bound]\label{theorem:DLRAMC_globalEB}
    Let $ \normalfont\textbf{Y}_{k}(\nu)$ denote the rank-$r$ approximation to $ \normalfont\boldsymbol{\Psi}_{k}(\nu)$ at $t_{k} = t_{0}+ kh$ obtained after $k$ steps of the augmented basis-update \& Galerkin integrator with step-size $h>0$ and $k=\frac{k^\prime}{h}$, with $k^\prime$ a fixed integer. Then, if the assumptions (A.1-A.3) are satisfied and $\vartheta_{i}(\nu)$ denotes the truncation error at $t_i$,  the error satisfies for all $k$ with $t_{k} =t_{0} + kh\leq t_{\mathrm{end}}$
    \begin{displaymath}
         \normalfont \lVert \mathbb{E}_{p}[\textbf{Y}_{k}(\nu) - \boldsymbol{\Psi}_{k}(\nu)] \rVert \leq \mathbb{E}_{p}[\widetilde{c}_{0}(\nu)] \delta + \mathbb{E}_{p}[\widetilde{c}_{1}(\nu)]\varepsilon + c_{2}h \,,
    \end{displaymath}
    where the constants $\widetilde{c_{0}}(\nu)$, $\widetilde{c_1}(\nu)$ and $c_2$ only depend on $L$, $B$, and $t_{\mathrm{end}}$. In particular, the constants are independent of singular values of the exact or approximate solution.
\end{theorem}

Thus, the error contribution to the bias due to low-rank approximation remains robust to small singular values in the probabilistic setting, making it a suitable and cheaper alternative to full-rank numerical solvers for uncertainty quantification. In Section~\ref{sec:DLRAMC_CVEst}, we discuss using reduced rank DLRA models as a control variate for the estimator $\widehat{Q}_{r;\mathrm{MC}}$.
\begin{remark}
    Note that the Lipshitz constant $L$ for the radiation transport equation~\eqref{eq:DLRAMC_3Dradiativetransport} is large, however, the above theorem gives insights into the behavior of the bias term. Additionally, when the chosen rank $r$ is greater than or equal to the true rank of the solution and the initial condition is also of rank $r$, the initial error $\delta = 0$.
\end{remark}

% The assumptions of the probabilistic robust error bound further imply that the low-rank approximation is always $\varepsilon$-close to the full (discretized) solution, independent of the random variable $\nu$. 

%%%%%%%%%%%%%%%%%%%%%%%%%%%%%%%%%%%%%%%%%%%%%%%%%%%%%%%%%%%%%%%%%%%%%%%%%%
\subsection{Control variates} \label{sec:DLRAMC_CVEst}
%%%%%%%%%%%%%%%%%%%%%%%%%%%%%%%%%%%%%%%%%%%%%%%%%%%%%%%%%%%%%%%%%%%%%%%%%%
The DLRA-based integrator reduces the computational cost of each Monte Carlo sample by using an approximate solution that is $\varepsilon$-close to the low-rank manifold independent of the random variable $\nu$. We thus expect to be able to compute more samples, at the same cost, compared to a full-rank integrator. Given that the variance of the MC estimator scales as $\mathcal{O}({N^{-1}})$, we expect to produce a lower variance estimate $\widehat{Q}_{r;\mathrm{MC}}$. However, a DLRA-based integrator capturing the true rank of the PDE solution may still be quite expensive in practice. This is further compounded by the fact that one often over-approximates the rank of the solution, out of caution. Hence, for high-dimensional PDEs, like \eqref{eq:DLRAMC_3Dradiativetransport}, it is infeasible to arbitrarily increase the number of MC samples. We therefore introduce a control-variate strategy to further suppress the variance of $\widehat{Q}_{r;\mathrm{MC}}$.

Simulations based on a reduced rank ($s<r$) can produce good, lower-cost approximations to more accurate higher-rank simulations. Hence, we can use the reduced-rank simulations as a control variate, with the goal of reducing the variance of the estimator $\widehat{Q}_{r;\mathrm{MC}}$. Let $\textbf{Y}^{r}_{k}(\nu)$ and $\textbf{Y}^{s}_{k}(\nu)$ be the low-rank approximations with ranks $r$ and $s$, respectively. We assume that $\mathbb{E}_{p}\left[\mathcal{G}(\textbf{Y}^{s}_{k}(\nu))\right]$ is known. Then for some $\alpha>0$, we define 
\begin{displaymath}
    \mathcal{G}_{\mathrm{CV}}(\nu) \coloneqq \mathcal{G}(\textbf{Y}_{k}^{r}(\nu)) - \alpha\left(\mathcal{G}(\textbf{Y}_{k}^{s}(\nu)) - \mathbb{E}_{p}\left[\mathcal{G}(\textbf{Y}_{k}^{s}(\nu))\right]\right).
\end{displaymath}
With re-arrangement we get the control variate estimator
\begin{equation}\label{eq:DLRAMC_CVestimator}
    \widehat{Q}_{r;\mathrm{CV}} = \alpha \mathbb{E}_{p}\left[\mathcal{G}(\textbf{Y}^{s}_{k}(\nu))\right] +  \frac{1}{N}\sum_{i=1}^{N}\left( \mathcal{G}(\textbf{Y}_{k}^{r}(\nu_{i})) - \alpha\mathcal{G}(\textbf{Y}_{k}^{s}(\nu_{i}))\right).
\end{equation}
If $\mathrm{Var}_{j}\coloneqq\mathrm{Var}(\mathcal{G}(\textbf{Y}_{k}^{j}(\nu)))$, $j\in\{r,s\}$, then the variance of $\mathcal{G}_{\mathrm{CV}}(\nu)$ is
\begin{displaymath}
        \mathrm{Var}(\mathcal{G}_{\mathrm{CV}}(\nu)) = \mathrm{Var}_{r} + \alpha^{2}\mathrm{Var}_{s}- 2\alpha\mathrm{Cov}_{rs},
\end{displaymath}
% \begin{displaymath}
%     \begin{split}
%         \mathrm{Var}(\mathcal{G}_{\mathrm{CV}}(\nu)) = \mathrm{Var}(\mathcal{G}(\textbf{Y}_{k}^{r}(\nu))) + \alpha^{2}\mathrm{Var}(\mathcal{G}(\textbf{Y}_{k}^{s}(\nu))) \\- 2\alpha\mathrm{Cov}(\mathcal{G}(\textbf{Y}_{k}^{r}(\nu)),\mathcal{G}(\textbf{Y}_{k}^{s}(\nu))),
%     \end{split}
% \end{displaymath}
\begin{displaymath}
\text{where} \qquad    \mathrm{Cov}_{rs} = \mathbb{E}_{p}\left[ \left\langle \mathcal{G}(\textbf{Y}_{k}^{r}(\nu)) - Q_{r},\mathcal{G}(\textbf{Y}_{k}^{s}(\nu)) - Q_{s} \right\rangle \right].
\end{displaymath}
% \begin{displaymath}
%     \mathrm{Cov}(\mathcal{G}(\textbf{Y}_{k}^{r}(\nu)),\mathcal{G}(\textbf{Y}_{k}^{s}(\nu))) = \mathbb{E}_{p}\left[ \left\langle \mathcal{G}(\textbf{Y}_{k}^{r}(\nu)) - Q_{r},\mathcal{G}(\textbf{Y}_{k}^{s}(\nu)) - Q_{s} \right\rangle \right].
% \end{displaymath}
The optimal $\alpha$, that minimizes the variance of the estimator is given by $ \alpha^{\ast} = \mathrm{Cov}_{rs}/\mathrm{Var}_{s}.$
\begin{remark}
    In practice, the optimal value of $\alpha$ as well as $\mathbb{E}_{p}\left[\mathcal{G}(\textbf{Y}^{s}_{k}(\nu))\right]$ have to be estimated (on-the-fly or beforehand). As $s<r$, the estimate is cheaper than the equivalent computation with rank-$r$. Furthermore, the covariance of $\mathcal{G}(\textbf{Y}_{k}^{r}(\nu))$ and $\mathcal{G}(\textbf{Y}_{k}^{s}(\nu))$ increases as $s$ approaches $r$.
\end{remark}

%%%%%%%%%%%%%%%%%%%%%%%%%%%%%%%%%%%%%%%%%%%%%%%%%%%%%%%%%%%%%%%%%%%%%%%%%%
\section{Numerical experiments} \label{sec:DLRAMC_numerical}
%%%%%%%%%%%%%%%%%%%%%%%%%%%%%%%%%%%%%%%%%%%%%%%%%%%%%%%%%%%%%%%%%%%%%%%%%%
To demonstrate the efficacy of the proposed low-rank estimators, we consider the 1D radiation transport equation under uncertain initial conditions. The source code can be found at \href{https://github.com/chinsp/publication-Low-rank-for-uncertain-radiative-transfer}{github.com/chinsp/publication-Low-rank-for-uncertain-radiative-transfer}. Specifically, we consider the radiation transport equation~\eqref{eq:DLRAMC_3Dradiativetransport} in slab geometry \cite{bell_nuclear_1970}, which assumes that the solution is rotationally invariant along infinite slabs perpendicular to the x-axis. This yields $ x\in [a,b] $ and the travel direction being restricted to $ \mu\in [-1,1] $. We further introduce an uncertain parameter $\nu$ with probability density function $p(\nu)$, giving the model equation
\begin{subequations}\label{eq:DLRAMC_RTslab}
    \begin{displaymath}
        \partial_{t}\psi(t,x,\mu,\nu) + \mu\partial_{x}\psi(t,x,\mu,\nu) = \sigma_{s}(x)\left( \frac{1}{2}\int_{-1}^{1}\psi(t,x,\mu',\nu) \dd\mu' - \psi(t,x,\mu,\nu) \right),
    \end{displaymath}
    \begin{displaymath}
         \psi(t_{0},x,\mu,\nu) = \psi_{0}(x,\mu,\nu).
    \end{displaymath}
\end{subequations}    
Note that we assume that the particle density never reaches the boundary. In the numerical experiments, the initial condition is set as a cut-off Gaussian with an uncertain amplitude:
\begin{displaymath}
    \psi(0,x,\mu,\nu) = \mathrm{max}\left\{ 10^{-4}, \frac{\nu}{\sqrt{2\pi}\sigma}\exp{\left(-\frac{x^{2}}{2\sigma^{2}}\right)} \right\},\quad x\in[-1.5,1.5],\! \nu\sim \mathrm{U}(0.5,1.5).
\end{displaymath}
The equation is discretized using $m$ grid points in space with an upwind-downwind flux and a spherical harmonics method of order $n-1$ (P$_{n-1}$) \cite{case_linear_1967} in the direction of travel. The discretization parameters are specified case-wise. This yields a matrix-valued differential equation with an uncertain initial condition. 

In this work, we are interested in the expected value of the scalar flux, defined as
\begin{displaymath}
    \mathcal{\phi}(x,\nu) \coloneqq \frac{1}{\sqrt{2}}\int_{-1}^{1}\psi(t=1.0,x,\mu,\nu)\hspace{1mm}\dd\mu,
\end{displaymath} 
i.e., $Q = \mathbb{E}_{p}\left[ \mathcal{\phi}(x,\nu) \right]$. A fine-grid Monte Carlo reference solution is computed for $\mathrm{P}_{101}$ approximation with $m = 1601$ spatial grid points and $N = 1.024\times 10^{5}$ samples. The CFL number, representing the number of grid cells traveled per time step, is set to $1$ in all experiments.

% \paragraph{Test case 1 (TC1)} 

% \paragraph{Test case 2 (TC2)} In the first test case the initial condition of the particle density is set as a cut-off Gaussian with uncertain center
% \begin{displaymath}
%     f(0,x,\mu,\nu) = \mathrm{floor}\left\{ 10^{-4}, \frac{1}{\sqrt{2\pi}\sigma}\exp{\left(-\frac{(x-\nu)^{2}}{2\sigma^{2}}\right)} \right\},
% \end{displaymath}
% where $\nu\in\mathrm{U}(-0.1,0.1)$.

%%%%%%%%%%%%%%%%%%%%%%%%%%%%%%%%%%%%%%%%%%%%%%%%%%%%%%%%%%%%%%%%%%%%%%%%%%
\subsection{Low-rank Monte Carlo}\label{sec:DLRAMC_numericalExperiments_MC}
%%%%%%%%%%%%%%%%%%%%%%%%%%%%%%%%%%%%%%%%%%%%%%%%%%%%%%%%%%%%%%%%%%%%%%%%%%

\begin{figure}[h]
    % \centering
    \begin{subfigure}[t]{0.3\linewidth}
        \includegraphics[width=\linewidth]{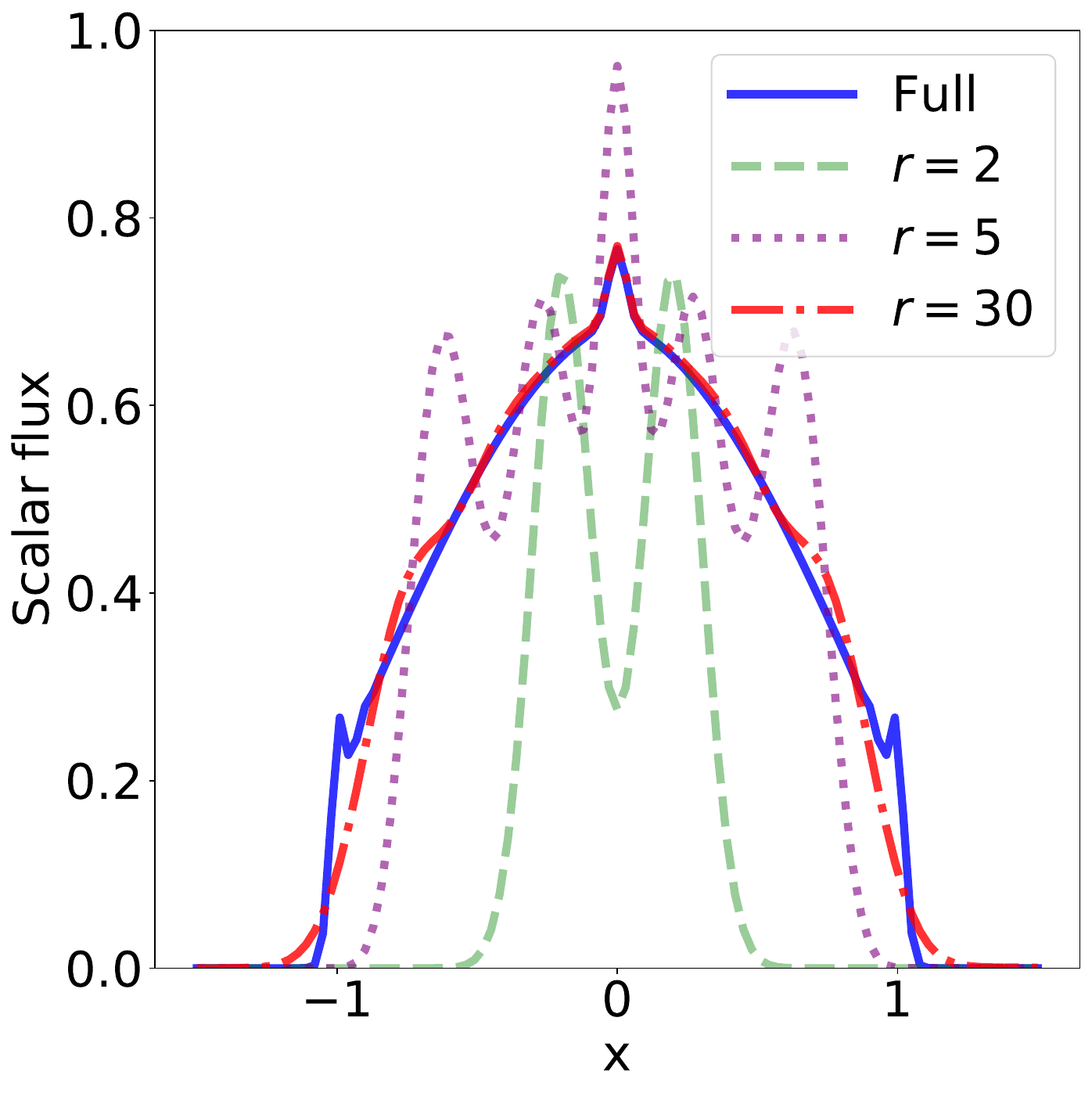}
        \caption{Low resolution grid}
        \label{fig:DLRAMC_Low_res_Scalar_flux}
    \end{subfigure}
    \begin{subfigure}[t]{0.3\linewidth}
        \includegraphics[width=\linewidth]{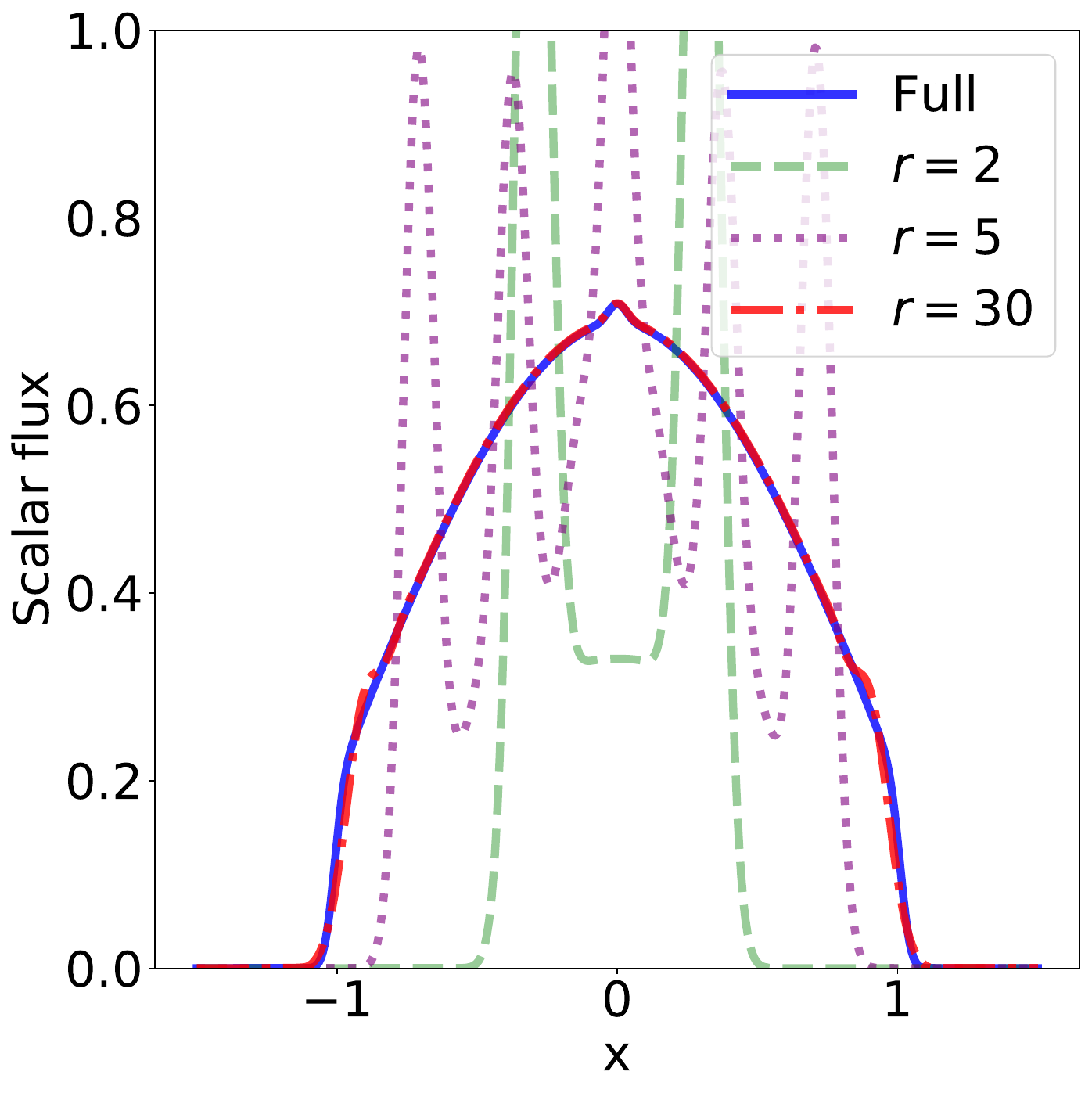}
        \caption{High resolution grid}
        \label{fig:DLRAMC_High_res_Scalar_flux}
    \end{subfigure}
    \begin{subfigure}[t]{0.3\linewidth}
      \centering
        \includegraphics[width=\linewidth]{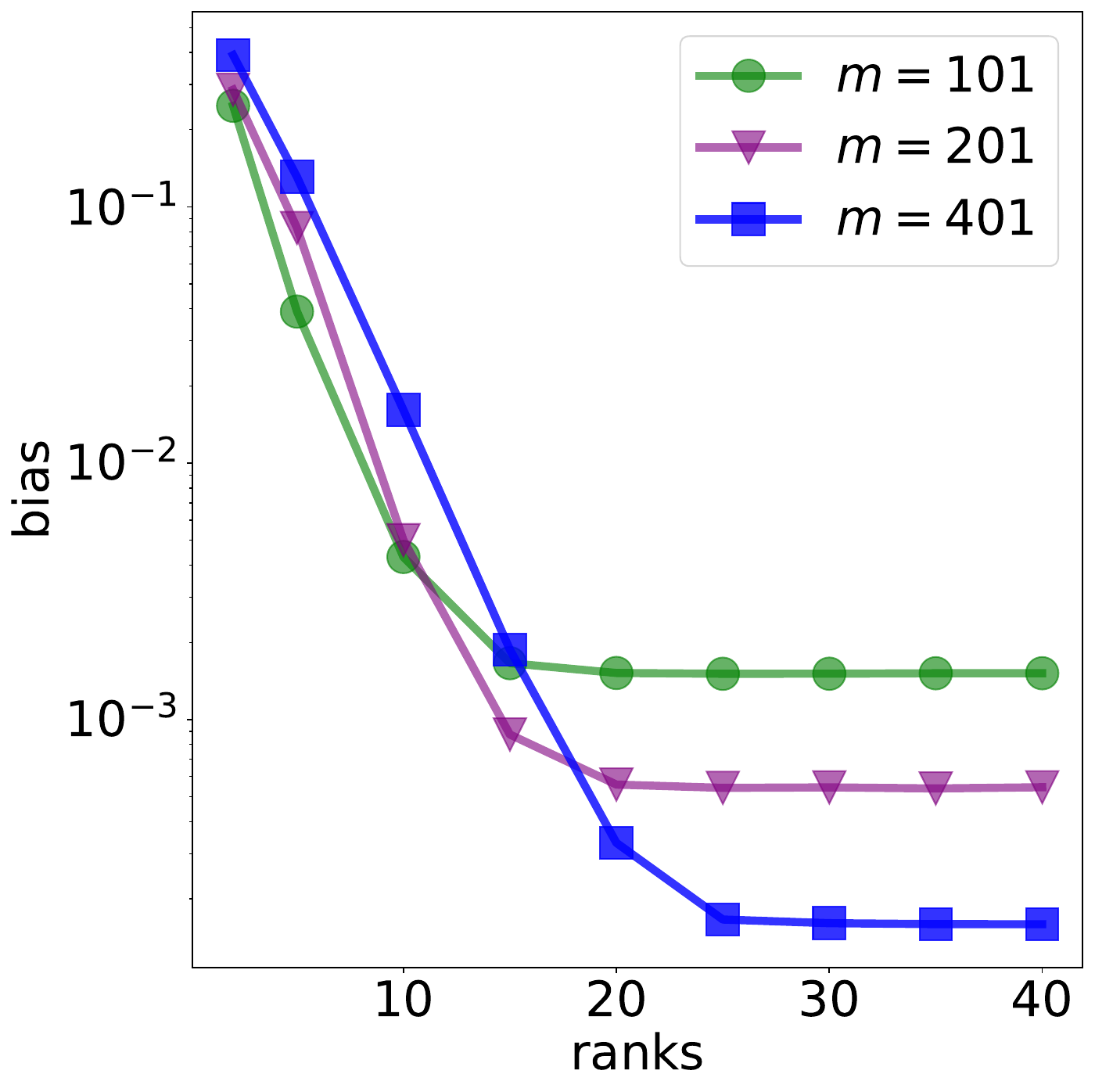}
        \caption{Bias}
        \label{fig:DLRAMC_bias}
    \end{subfigure}
    
    \begin{subfigure}[t]{0.3\textwidth}
        \includegraphics[width=\textwidth]{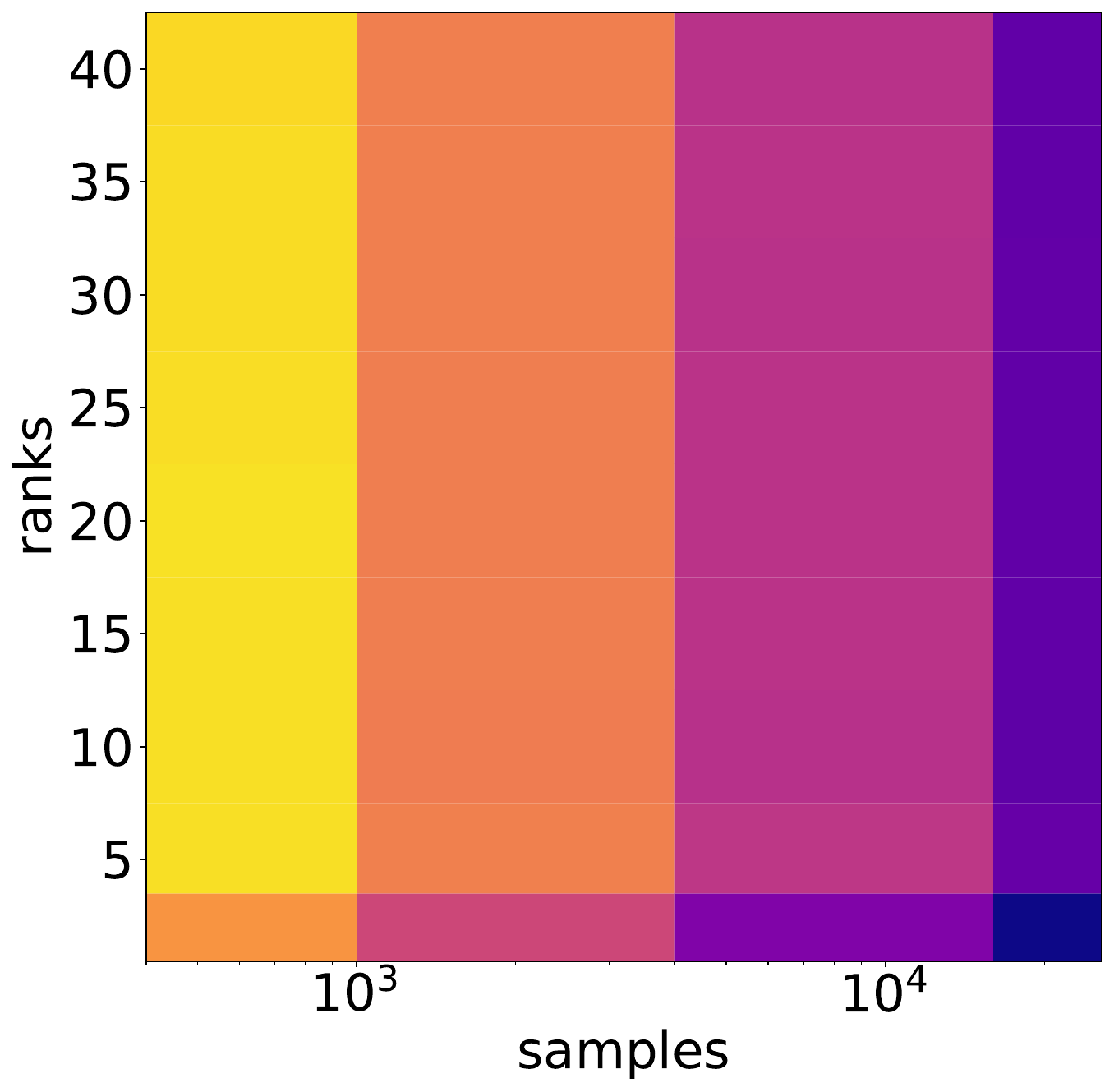}
        \caption{$m = 101$}
        \label{fig:DLRAMC_MCerror_m100}
    \end{subfigure}
    \begin{subfigure}[t]{0.3\textwidth}
        \includegraphics[width=\textwidth]{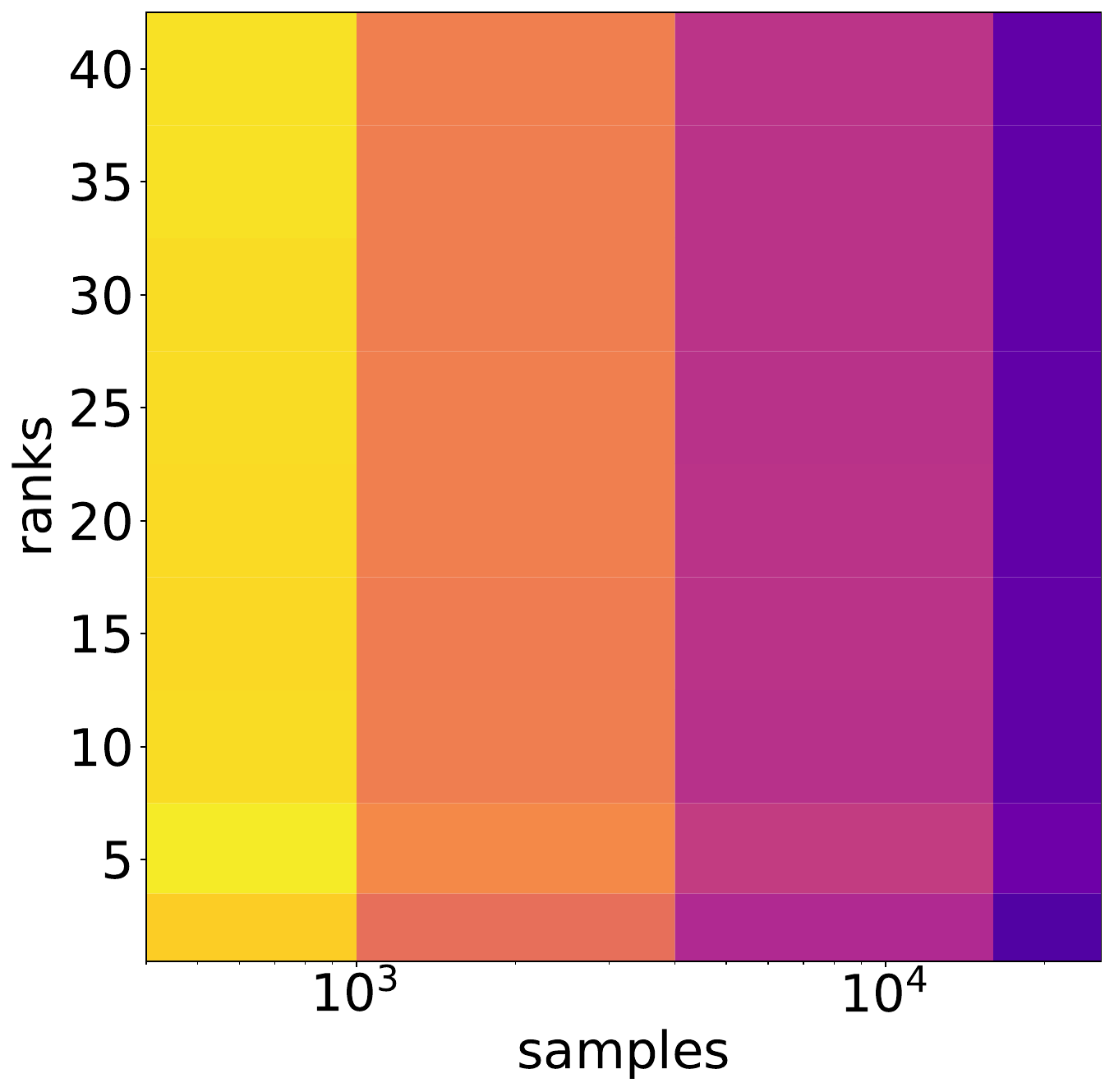}
        \caption{$m = 201$}
        \label{fig:DLRAMC_MCerror_m200}
    \end{subfigure}
    \begin{subfigure}[t]{0.3\textwidth}
        \includegraphics[width=\textwidth]{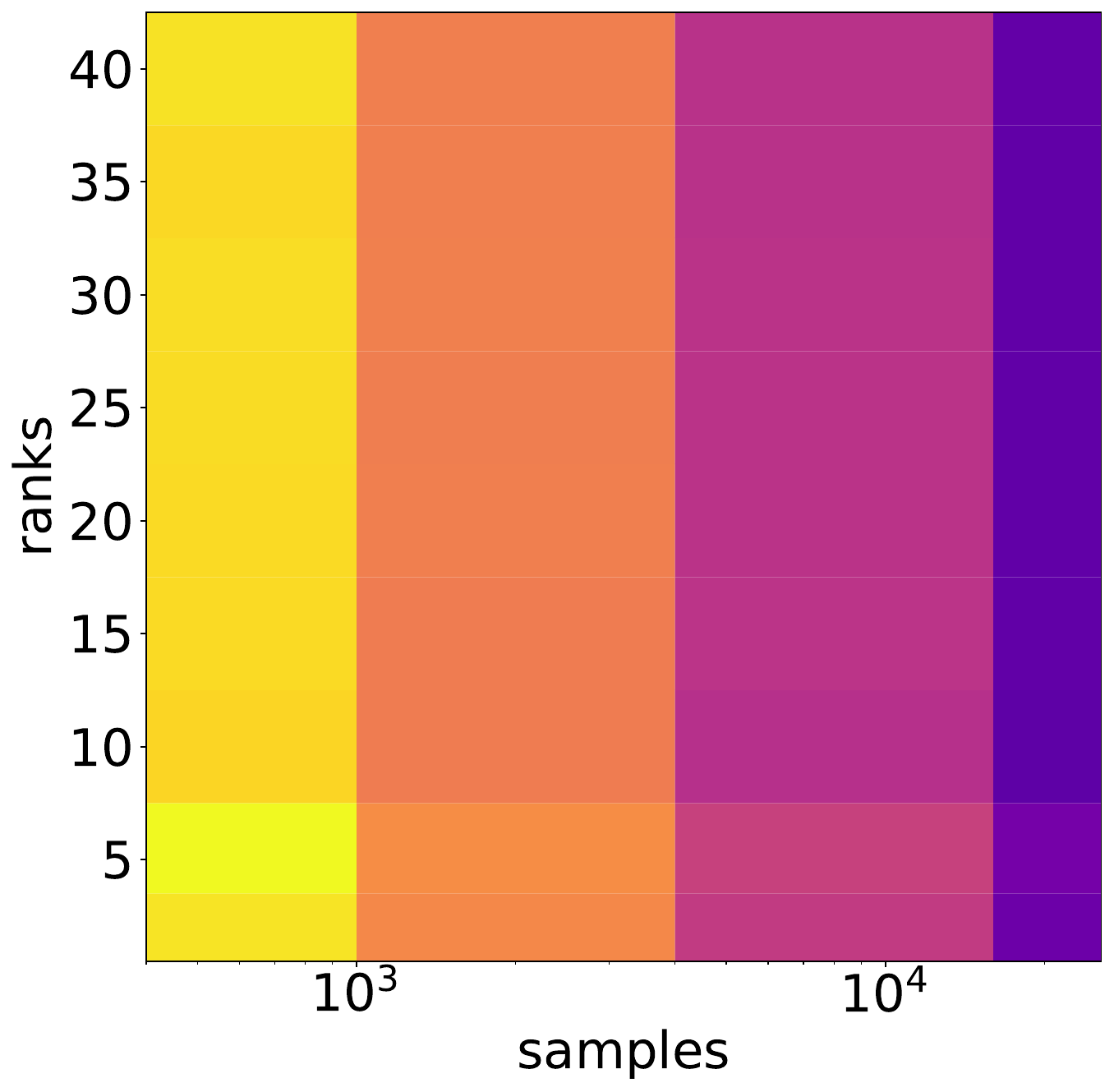}
        \caption{$m = 401$}
        \label{fig:DLRAMC_MCerror_m400}
    \end{subfigure}
    \begin{subfigure}[t]{0.08\textwidth}
    \vspace*{-3.6cm}\includegraphics[width=\linewidth,height=3.23cm,keepaspectratio]{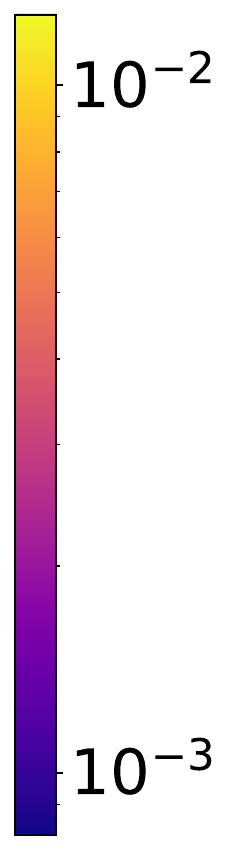}
    \end{subfigure}
    
   % \begin{subfigure}[t]{0.3\textwidth}
   %      \includegraphics[width=\textwidth]{figures/BUGvsFull_rank_samples_bias_normalized_101.pdf}
   %      \caption{$m = 100$}
   %  \end{subfigure}
   %  \begin{subfigure}[t]{0.3\textwidth}
   %      \includegraphics[width=\textwidth]{figures/BUGvsFull_rank_samples_bias_normalized_201.pdf}
   %              \caption{$m = 200$}
   %  \end{subfigure}
   %  \begin{subfigure}[t]{0.3\textwidth}
   %      \includegraphics[width=\textwidth]{figures/BUGvsFull_rank_samples_bias_normalized_401.pdf}
   %      \caption{$m = 400$}
   %  \end{subfigure}
   %  \begin{subfigure}[t]{0.08\linewidth}     \vspace*{-3.68cm}\includegraphics[width=\linewidth,height=3.32cm,keepaspectratio]{figures/colorbar_bias.pdf}
   %  \end{subfigure}

    \caption{Top row: Scalar flux computed using the full P$_{101}$ solver and the augment BUG solver at time $t=1$ for a fixed sample $\nu$ and ranks $r\in\{2,5,30\}$ with (a) $m=101$ (b) $m=401$ spatial grid points. (c) The bias of the low-rank Monte Carlo estimator for $r=30$ and $N=25600$ samples. Bottom row: Monte Carlo error of the low-rank Monte Carlo estimator in log-scale for spatial grids of varying widths. }
    \label{fig:DLRAMC_MCerror_bias}
\end{figure}

While applying the low-rank Monte Carlo estimator is straightforward, it is unclear how the estimator's parameters affect the MC error and the bias of the estimate. Thus, we conduct a parametric study of the low-rank Monte Carlo method by varying the spatial grids, rank, and the Monte Carlo samples. In preliminary investigations, it was seen that the discretization of $\mu$ does not affect the Monte Carlo error or the bias and thus we fix the discretization to the same one as the reference. For the spatial discretization we use nested grids with $m \in \{101,201,401\}$, ranks $r \in \{2,5,10,15,20,25,30,35,40\}$, and $N\in\{ 400,1600,6400,25600 \}$, where $N$ is the number of samples. The results are plotted in Figure~\ref{fig:DLRAMC_MCerror_bias}.

We see from Figure~\ref{fig:DLRAMC_Low_res_Scalar_flux}-\ref{fig:DLRAMC_bias} that, independent of the underlying spatial grid, the bias decreases with increasing rank until the spatial discretization error dominates it. Additionally, we see that for smaller ranks ($r\leq 10$), the bias is smaller when the underlying spatial grid is coarse. At the same time, the behavior is flipped for larger ranks ($r\geq 20$), i.e., the error is smaller when the underlying spatial grid is finer. From Figure~\ref{fig:DLRAMC_MCerror_m100}-\ref{fig:DLRAMC_MCerror_m400} we see that as the sample size increases, the Monte Carlo error decreases independent of the rank or the underlying spatial grid. Hence, a good estimate of the QoI requires knowledge of the interplay of rank and the underlying spatial discretization. Here, especially the required rank is difficult to predict in practice, often resulting in over-approximation. Further, even though each individual sample is cheaper to compute using a low-rank approach, reducing the error of the low-rank Monte Carlo estimator can still require a large number of samples. Especially for higher-dimensional problems like \eqref{eq:DLRAMC_3Dradiativetransport} it is therefore worthwhile to consider further approaches to variance reduction within our low-rank framework.

%%%%%%%%%%%%%%%%%%%%%%%%%%%%%%%%%%%%%%%%%%%%%%%%%%%%%%%%%%%%%%%%%%%%%%%%%%
\subsection{Control variates}
%%%%%%%%%%%%%%%%%%%%%%%%%%%%%%%%%%%%%%%%%%%%%%%%%%%%%%%%%%%%%%%%%%%%%%%%%%
%%%%%%%%%%%%%%%%%%%%%%%%%%%%%%%%%%%%%%%%%%%%%%%%%%%%%%%%%%%%%%%%%%%%%%%%%%
To gauge the potential benefits of using lower-rank estimates for variance reduction, we next compare the low-rank Monte Carlo approach to control variate estimators with different levels of coarseness with respect to rank. For this, we perform low-rank Monte Carlo samples as described in Section~\ref{sec:DLRAMC_numericalExperiments_MC}, with a fixed discretization of $x$ and $\mu$ ($m=201, n=101$) but varying rank $r \in \{20,25,30,35,40\}$. The number of samples is also fixed to $N_{\mathrm{MC}}=2000$. To improve the variance of the low-rank Monte Carlo estimator, we use a control variate with a coarser rank $s\in \{2,5,10,15\}$. Since the ultimate goal of variance reduction is typically to achieve a lower error at the same computational costs or vice versa, we fix the error and investigate the effect on the runtime. For this, we use the following heuristic: We compute the same number of samples for the coarse rank $s$ as for the low-rank Monte Carlo estimator, i.e., $N_{\mathrm{CV}_{\mathrm{coarse}}} = N_{\mathrm{MC}} = 2000$. Then, assuming the optimal $\alpha^{\ast}$ is known exactly, we can compute the optimal number of samples for the differences such that the total error of the control variate estimator is at most as high as for the rank-$r$ Monte Carlo estimator:
\begin{displaymath}
    N_{\mathrm{CV}_{\mathrm{diff}}} = \mathrm{Var}_{r}\cdot\frac{1-\mathrm{Corr}_{rs}^{2}}{\epsilon^{2}}, \quad \mathrm{Corr}_{rs} = \frac{\mathrm{Cov}_{rs}}{\sqrt{\mathrm{Var}_{r}}\sqrt{\mathrm{Var}_{s}}},
\end{displaymath}
where $\epsilon$ is the Monte Carlo error of $N_{\mathrm{MC}}$ samples (see~\ref{sec:DLRAMC_numericalExperiments_MC}), $\mathrm{Var}_{r}$ is the variance of the rank-$r$ solver, $\mathrm{Corr}_{rs}$ is  correlation and  $\mathrm{Cov}_{rs}$ the covariance between the finer rank-$r$ and coarser rank-$s$ solver. Since the above formula only holds when we know $\alpha^{\ast}$ exactly, we gauge the influence of having a priori knowledge of $\alpha^{\ast}$ on the estimated $N_{\mathrm{CV}_{\mathrm{diff}}}$ and thereby on the runtime of the control variate estimator. We do so by considering the following two approaches: in the first approach, we obtain an a priori estimate of $\alpha^{\ast}$ by running 500 MC samples for the coarse and fine rank approximations. This optimal value is then used to compute the optimal number of samples for the control variate estimator. In the second approach, we estimate $\alpha^{\ast}$ during runtime by using $N_{\mathrm{w}}$ warm-up samples based on which we estimate $N_{\mathrm{CV}_{\mathrm{diff}}}$. If $N_{\mathrm{CV}_{\mathrm{diff}}} \leq N_{\mathrm{w}}$, we re-use the estimates obtained from the warm-up samples and only compute the remaining $(N_{\mathrm{CV}} - N_{\mathrm{w}})$ samples for the coarse level. On the other hand, if $N_{\mathrm{CV}_{\mathrm{diff}}} > N_{\mathrm{w}}$, we compute $(N_{\mathrm{CV}_{\mathrm{diff}}} - N_{\mathrm{w}})$ extra samples for the difference $\left( \mathcal{G}(\textbf{Y}_{k}^{r}(\nu)) - \alpha^{\ast}\mathcal{G}(\textbf{Y}_{k}^{s}(\nu))\right)$ and $N_{\mathrm{CV}} - (N_{\mathrm{CV}_{\mathrm{diff}}} - N_{\mathrm{w}})$ for the coarse level. Note that the second approach often leads to a higher than necessary number of samples on the differences if we use $N_{\mathrm{w}}< 200$ warm-up samples and thus we set $N_{\mathrm{w}} = 200$ in all the experiments.  
\begin{figure}[h]
    \centering
    \begin{subfigure}{0.475\textwidth}
        \includegraphics[width=\linewidth]{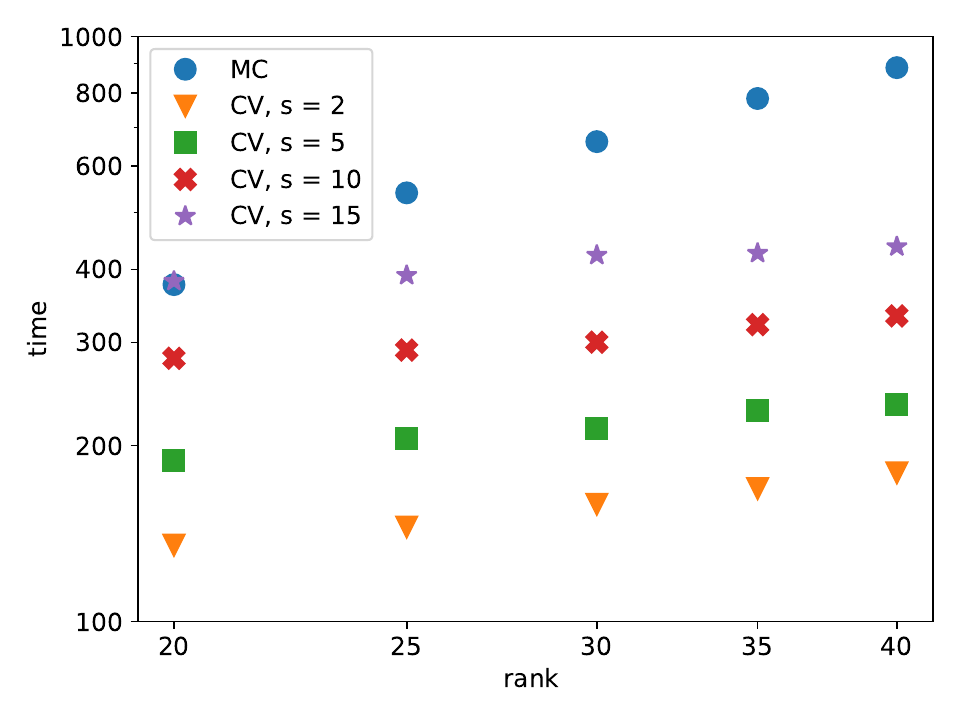}
        \caption{}
        \label{fig:timeCVMC_a}
    \end{subfigure}
      \begin{subfigure}{0.475\textwidth}
           \includegraphics[width=\linewidth]{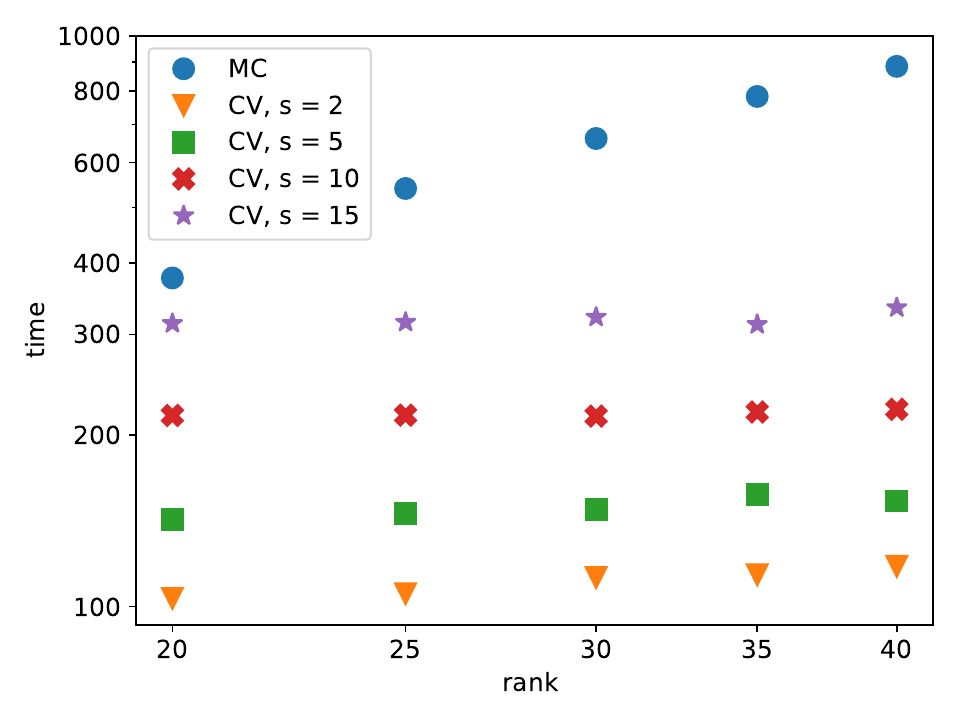}
        \caption{}
        \label{fig:timeCVMC_b}
    \end{subfigure}
    \caption{Log-log plot of (minimal) runtime of Monte Carlo vs. different control variates at theoretically constant error. (a) Number of samples and $\alpha$ were estimated based on a separate simulation with 500 samples. (b) Number of samples and $\alpha$ were estimated during runtime based on 200 warm-up samples.}
    \label{fig:timeCVMC}
\end{figure}
In Figure~\ref{fig:timeCVMC}, the minimal runtime over 20 runs is plotted for both the Monte Carlo estimator at different ranks and the control variate estimators with all combinations of coarse and fine ranks. Note that the runtime for the control variate also includes time taken to compute $N_{\mathrm{CV}}$ samples of the rank-$s$ approximation similar to a MLMC approach. It is apparent that the control variate estimators are consistently faster than low-rank Monte Carlo, both when computing the optimal number of samples beforehand (Figure~\ref{fig:timeCVMC_a}) and when estimating the parameters during computation using 200 warm-up samples (Figure~\ref{fig:timeCVMC_b}). While the runtime seems to increase exponentially with the growing rank for Monte Carlo, the runtimes for control variates grow much slower, closer to linearly, especially for the higher choices of coarse rank. This implies that it might be worthwhile to consider a multi-level or at least multiple level framework. 

Lastly, in Table~\ref{tab:optAlpha}, the precomputed optimal values for the parameter $\alpha$ are shown. While the values are almost constant with respect to the rank of the finer level $r$, they vary between around 0.6 and close to 1 depending on the rank of the coarse level. As expected, the correlation between fine and coarse ranks is higher for larger ranks on the coarse level. Thus, the $\alpha^{\ast}$ values are lower for smaller ranks, especially for $s = 2$ and $s = 4$ and approach 1 for the higher ranks.
\begin{table}[h]
    \centering
     \begin{tabular}{cccccc}
       $s$ $\setminus$ \raisebox{0.1cm}{$r$} &20 &25&30&35&40  \\
        \hline 
       2 $\quad$ &  0.6060 & 0.6060 & 0.6060 & 0.6060 & 0.6060\\
        5 $\quad$ & 0.7834 & 0.7833 & 0.7833 & 0.7833 & 0.7833\\
        10 $\quad$ & 0.9937 & 0.9940 & 0.9940 & 0.9940 & 0.9940\\
        15 $\quad$ & 0.9979 & 0.9980 & 0.9980 &  0.9980 &0.9980\\
         \hline \\
    \end{tabular}
    \caption{Precomputed $\alpha^{\ast}$ values obtained from the first approach with 500 samples.}
    \label{tab:optAlpha}
\end{table}
%\begin{figure}[h]
 %   \centering
 %       \includegraphics[width=0.75\linewidth]{figures/optAlpha_BoxPlot.png}

%    \caption{Violin plot of chosen optimal $\alpha$ over 20 runs for different ranks. If $\alpha>1$ it was set to 1 for the computation.}
 %%\end{figure}

\section{Conclusion} \label{sec:DLRAMC_conclusion}
We have studied the use of the dynamical low-rank approximation in the context of uncertainty quantification, both for cost reduction in a standard Monte Carlo framework and as a lower fidelity control variate. We showed that the robust error bound holds in a probabilistic setting and from the low-rank Monte Carlo experiments we see that, while the variance is only influenced by the sample size, an optimal balance of rank and grid size is necessary to correctly capture the quantity of interest. Computational costs can be further reduced significantly compared to the low-rank Monte Carlo method using a control variate approach based on a lower rank approximation. Based on our results, future work on the design of a fully multi-level approach balancing both rank and spatial discretization error seems promising for tackling larger scale problems involving the radiative transfer equation, as well as other problems with an identified low-rank structure.
%%%%%%%%%%%%%%%%%%%%%%%%%%%%%%%%%%%%%%%%%%%%%%%%%%%%%%%%%%%%%%%%%%%%%%%%%%

%%% To ensure the bibliography has the correct style please run bibtex with the spmpsci style
%%% which is included in the Springer zip file
%%% Here we assume refs.bib would be the name of the bib file containing bibliographic info
%%% You can then copy the .bbl produced, as given in the example below
%%%
\section*{Acknowledgements}
The work of Chinmay Patwardhan was funded by the Deutsche Forschungsgemeinschaft (DFG, German Research Foundation) – Project-ID 258734477 – SFB 1173. Pia Stammer received funding from the German National Academy of Sciences Leopoldina for the project underlying this article, under grant number LPDS 2024-03. 

\bibliographystyle{styles/bibtex/spmpsci}
%\bibliography{refs}
% \bibliographystyle{plain}
\bibliography{references}

\end{document}